\documentclass[12pt,leqno]{article}

\usepackage{amsthm}
\usepackage[usenames]{color}

\usepackage{amsmath}
\usepackage{amscd}
\usepackage{amsmath}
\usepackage{amssymb}
\usepackage{latexsym}
\makeatletter

\@addtoreset{equation}{section}
\makeatother

\newtheorem{thm}{Theorem}[subsection]
\newtheorem{pr}[thm]{Proposition}
\newtheorem{df}[thm]{Definition}
\newtheorem{lm}[thm]{Lemma}
\newtheorem{cor}[thm]{Corollary}
\newtheorem{cn}[thm]{Conjecture}

\newtheorem{ex}[thm]{Example}
\newtheorem*{cl}{Claim}

\newcommand{\sm}{\raisebox{2.33pt}{~\rule{6.4pt}{1.3pt}~}}

\input xy \xyoption{all} \CompileMatrices

\begin{document}

\title{Characteristic cycles
and the conductor 
of direct image}
\author{Takeshi Saito}

\maketitle

\begin{abstract}
We prove the functoriality
for proper push-forward
of the characteristic cycles 
of constructible complexes
by morphisms of
smooth projective schemes over
a perfect field,
under the assumption that 
the direct image of the
singular support has the
dimension at most
that of the target
of the morphism.
The functoriality 
is deduced from a conductor
formula which is a special case
for morphisms to curves.
The conductor formula in the
constant coefficient case gives
the geometric case of
a formula conjectured by Bloch.
\end{abstract}

Let $k$ be a perfect field
and $\Lambda$ be a finite field
of characteristic invertible in $k$.
For a constructible complex
${\cal F}$ of $\Lambda$-modules
on a smooth scheme $X$ over $k$,
the characteristic cycle
$CC{\cal F}$ is defined in \cite[Definition 5.10]{CC}
as a cycle supported on
the singular support $SS{\cal F}$
defined by Beilinson in \cite{Be}
as a closed conical subset of the cotangent bundle $T^*X$.
We study the functoriality
of characteristic cycles 
for proper push-forward.

Let $f\colon X\to Y$
be a morphism
of smooth projective 
schemes over $k$.
Then, we prove in
Theorem \ref{thmf*}
the equality
\begin{equation*}
CCRf_*{\cal F}=
f_!CC{\cal F}
\leqno{\rm (\ref{eqcnf})}
\end{equation*}
conjectured in \cite[Conjecture 1]{prop}
under the assumption
$\dim f_\circ SS{\cal F}\leqq \dim Y=m$
for the direct image
$f_\circ SS{\cal F}\subset T^*Y$.
The precise definitions
will be given in Subsection \ref{ssdc}.
We can slightly weaken the assumption, as is seen in Theorem
\ref{thmf*}.

A typical example where the assumption 
$\dim f_\circ SS{\cal F}\leqq \dim Y=m$
fails is the case where $f$ is the Frobenius.
Without the assumption, the author does not know
how to prove the equality (\ref{eqcnf}) 
because the Milnor formula (\ref{eqMil}) on which
the definition of characteristic cycles is based
will not make sense.
In the general case,
the equality (\ref{eqcnf}) 
implies an equality in the Chow group
${\rm CH}_m(T^*Y)={\rm CH}_0(Y)$.
In the case where $Y={\rm Spec}\, k$ for a finite field $k$,
the equality in ${\rm CH}_0(Y)={\mathbf Z}$
is proved in \cite{UYZ} under the assumption
that $X$ is projective by a method different from
that in the present article using the product formula
for the constant term of the functional equation of
$L$-function.

The formula (\ref{eqcnf})
is an algebraic analogue of \cite[Proposition 9.4.2]{KSc} where 
functorial properties of characteristic cycles are studied in a transcendental context. 
In the case where $Y={\rm Spec}\ k$,
the equality (\ref{eqcnf}) is 
the index formula 
\begin{equation*}
\chi(X_{\bar k},{\cal F})
= (CC{\cal F},T^*_XX)_{T^*X}
\leqno{\rm (\ref{eqind})}
\end{equation*}
computing the Euler-Poincar\'e characteristic  as an intersection number
proved in
\cite[Theorem 7.13]{CC}.

We deduce the functoriality
(\ref{eqcnf})
from the index formula 
(\ref{eqind}) 
in Subsection \ref{ssccd} as follows.
By taking a projective embedding
of $Y$ and a good pencil,
we reduce it to the case
where $Y$ is
a projective smooth curve.
By the index formula
(\ref{eqind}) applied to a general fiber,
the equality (\ref{eqcnf})
is equivalent to a conductor formula 
\begin{equation*}
-a_yRf_*{\cal F}=
(CC{\cal F},df)_{T^*X,X_y}
\leqno{\rm (\ref{eqcf})}
\end{equation*}
proved in
Theorem \ref{thmay},
where the left hand side
denotes the Artin conductor
at a closed point $y\in Y$
of the direct image.
In the case where ${\cal F}$
is the constant sheaf $\Lambda$,
the right hand side equals
the localized self-intersection product
defined in \cite{Bl}
and the formula (\ref{eqcf})
specializes to the geometric case, 
Corollary \ref{corBl}, of
the conductor formula conjectured
in \cite{Bl} by Bloch.
In \cite{TV}, the authors announce an attempt
to prove the conductor formula
by a different approach.

Further the index formula implies
that we have an equality (\ref{eqsum})
for the sums over $y\in Y$
of the both sides
in (\ref{eqcf}).
To deduce (\ref{eqcf}) from (\ref{eqsum}) 
for the sums,
it suffices to show the existence of
a covering of $Y$
\'etale at a fixed point $y$
killing the contributions of the
other points.

For the vanishing of
the left hand side,
the local acyclicity of 
$f\colon X\to Y$ relative to ${\cal F}$ 
is a sufficient condition.
The $SS{\cal F}$-transversality of 
$f\colon X\to Y$ 
defined in Definition \ref{dftrans}
and studied in Subsection \ref{ssFtr}
after some preliminaries in
Subsection \ref{ssCtr}
is a stronger condition
and is a sufficient
condition for the vanishing of
the right hand side.
Thus, the proof of (\ref{eqcf}) 
is reduced to showing
variants of the stable reduction theorem
on the existence of ramified covering $Y'$ of $Y$
and of a modification ${\cal F}'$ 
of the pull-back of ${\cal F}$ on
$X'=X\times_YY'$ such that
the base change 
$f'\colon X'\to Y'$ is
locally acyclic relatively to ${\cal F}'$ 
and is $SS{\cal F}'$-transversal.

We show that the existence of 
a modification of a perverse sheaf
${\cal F}$ relatively to which 
$f\colon X\to Y$ is locally acyclic 
is equivalent to the condition that
the inertia action on
the nearby cycles complex
$R\Psi{\cal F}$ is trivial
in Proposition \ref{prS}.2.
This is rather a direct consequence
of the relation of
the direct image by the open 
immersion of the generic fiber
with the nearby cycles complex.
As we work with torsion coefficients,
the condition is satisfied over
a ramified covering of $Y$.

Further,
we show that
the local acyclicity of
$f\colon X\to Y$ relatively to ${\cal F}$
implies
the existence of a ramified covering $Y'\to Y$
such that the
base change of
$f\colon X\to Y$ is 
$SS{\cal F}'$-transversal 
for the pull-back ${\cal F}'$
of ${\cal F}$
in Corollary \ref{corV} 
of Theorem \ref{thmV2}.
Theorem \ref{thmV2} is deduced 
from a weaker version
Proposition \ref{prV2}
which is proved by using
the alteration \cite[Theorem 8.2]{dJ}.
In Proposition \ref{prV2},
the ramified covering may be inseparable,
while it is generically \'etale
in Theorem \ref{thmV2}.
This improvement is crucial because in the proof
of Theorem \ref{thmay},
we need to find a covering of a curve $Y$
which is \'etale at a fixed point $y\in Y$.
Theorem \ref{thmV2}
is proved by an argument similar to
that in the proof of 
\cite[Proposition 3.2]{TF}
by using a consequence of
the stable reduction theorem
\cite[Theorem 1.5]{Tem}.

We also prove an index formula
Proposition \ref{prvan}
for vanishing cycles complex.

The author thanks A.\ Beilinson
for the remark that 
Theorem \ref{thmay}
implies the geometric case
of the conductor formula
conjectured by Bloch in \cite{Bl}
and for showing
the proof of Lemma \ref{lmch0}
in the characteristic zero case.
The author thanks H.\ Haoyu
for discussion on the subject
of Subsection \ref{ssiv}
and thanks H.\ Kato
for pointing out an error
in the proof of Proposition \ref{prvan}
in an earlier version.
He also thanks an anonymous referee
for careful reading,
for pointing out gaps
in the proofs of Lemma \ref{lmtrqf}
and of Proposition \ref{prV2}
in an earlier version
and proposing to improve the statement of
Lemma \ref{lmS}.

The research was supported
by JSPS Grants-in-Aid 
for Scientific Research
(A) 26247002.

\tableofcontents

\section{Local acyclicity and transversality}\label{str}

\subsection{Preliminaries on perverse sheaves}
\label{ssperv}

We fix some conventions on
perverse sheaves.
Let $X$ be a noetherian
scheme and let $\Lambda$ be
a finite field of characteristic $\ell$
invertible on $X$.
We say that a complex ${\cal F}$ of
$\Lambda$-modules
on the \'etale site of $X$
is constructible
if the cohomology sheaf
${\cal H}^q{\cal F}$ is
constructible for every integer $q$
and 
${\cal H}^q{\cal F}=0$
except for finitely many $q$.
Let $D^b_c(X,\Lambda)$
denote the category of constructible
complexes of $\Lambda$-modules.

First we recall the case
where $X$ is a scheme of finite
type over a field $k$.
Let $\Lambda$ be
a finite field of characteristic $\ell$
invertible in $k$.
Then, the $t$-structure
$(^{\rm p}\!D^{\leqq 0},\,
^{\rm p}\!D^{\geqq 0})$
on $D^b_c(X,\Lambda)$
relative to the middle perversity
is defined in \cite[2.2.10]{BBD}
and the perverse sheaves
on $X$ form an abelian
subcategory
${\rm Perv}(X,\Lambda)
=\,
^{\rm p}\!D^{\leqq0}
\cap\,
^{\rm p}\!D^{\geqq0}.$

Every object of
${\rm Perv}(X,\Lambda)$
is of finite length by
\cite[Th\'eor\`eme 4.3.1 (i)]{BBD}.
Further by 
\cite[Th\'eor\`eme 4.3.1 (ii)]{BBD}, 
simple objects are of the following form:
Let $V\subset X$ be an irreducible 
locally closed subset of dimension $d$
such that
the reduced part of the geometric fiber
$V_{\bar k}$ is smooth.
Let $j\colon V\to X$ be the immersion
and let ${\cal G}$ be an irreducible
locally constant sheaf on $V$.
Then, $j_{!*}{\cal G}[d]$ is a simple perverse sheaf
on $X$.

Let $j\colon U\to X$ be an open immersion
and ${\cal F}$ be a perverse sheaf on $X$.
Let ${\cal G}\subset {\cal F}$
be the largest sub perverse sheaf
supported on the complement
$X\sm U$
and let
${\cal G}\subset {\cal H}\subset {\cal F}$
be the smallest sub perverse sheaf
such that ${\cal F}/{\cal H}$ 
is supported on $X\sm U$.
Then, since $j_{!*}j^*{\cal F}$ is
the unique extension of $j^*{\cal F}$
without non-trivial sub or quotient perverse sheaf
supported on $X\sm U$
by \cite[Corollaire 1.4.25]{BBD},
there is a canonical isomorphism
\begin{equation}
{\cal H}/{\cal G}\to j_{!*}j^*{\cal F}.
\label{eqj!*}
\end{equation}

\begin{pr}\label{prffp}
Let $h\colon W\to X$ be a finite
and faithfully flat morphism
of schemes of finite type over a field $k$
and let ${\cal F}\in
{\rm Perv}(X,\Lambda)$ be a perverse sheaf
on $X$.
Then, the trace morphism
$h_*h^*{\cal F}\to {\cal F}$
{\rm (\cite[Th\'eor\`eme 6.2.3]{SGA4})}
induces a surjection
\begin{equation}
h_*\, ^{\rm p}\!H^0 h^*{\cal F}\to {\cal F}
\label{eqTr}
\end{equation}
of perverse sheaves on $X$.
\end{pr}

\proof{
Since $h^*$ is right $t$-exact by
\cite[Proposition 2.2.5]{BBD}, 
and $h_*$ is $t$-exact by
\cite[Corollaire 4.1.3]{BBD},
for perverse sheaf ${\cal F}$ on $X$,
we have 
$h_*h^*{\cal F}\in D^{\leqq 0}$
and 
$^{\rm p}\!H^0h_*h^*{\cal F}=
h_*{}^{\rm p}\!H^0h^*{\cal F}$
and the trace morphism
induces (\ref{eqTr}).

First we show the surjectivity 
for a simple perverse sheaf
${\cal F}$ on $X$.
Since ${\cal F}$ is simple,
it suffices show that
(\ref{eqTr}) is non trivial.
By \cite[Th\'eor\`eme 4.3.1 (ii)]{BBD},
we may assume that
${\cal F}=j_{!*}{\cal G}[d]$,
for an irreducible 
locally closed subset $V\subset X$ of dimension $d$
such that
the reduced part of the geometric fiber
$V_{\bar k}$ is smooth,
the immersion
$j\colon V\to X$
and an irreducible
locally constant sheaf ${\cal G}$ on $V$.
Let $h_V\colon V'\to V$ denote the base change of $h$.
Shrinking $X$, we may assume
that $V\subset X$ is a closed subset.
Further shrinking $X$, we may assume that
the reduced part of $V'_{\bar k}$ is smooth
over $\bar k$ and
that $h_{V*}h_V^*{\cal G}$ is a locally
constant sheaf on $V$.
Then
$h^*{\cal F}$ is a perverse sheaf on $W$
and 
the trace morphism
$h_{V*}h_V^*{\cal G}\to {\cal G}$
is a surjection. Hence the assertion follows
in this case.

We show the general case.
Since $h^*$ is right $t$-exact by
\cite[Proposition 2.2.5]{BBD}, 
and $h_*$ is $t$-exact by
\cite[Corollaire 4.1.3]{BBD}, for an exact sequence
$0\to {\cal F}'\to {\cal F}\to {\cal F}''\to 0$
of perverse sheaves on $X$,
we have a commutative diagram 
\begin{equation}
\begin{CD}
@.h_*\, ^{\rm p}\!H^0 h^*{\cal F}'@>>> 
h_*\, ^{\rm p}\!H^0 h^*{\cal F}@>>> 
h_*\, ^{\rm p}\!H^0 h^*{\cal F}''@>>> 0\\
@.@VVV@VVV@VVV@.\\
0@>>> {\cal F}'@>>> {\cal F}@>>> {\cal F}''@>>> 0
\end{CD}
\end{equation}
of exact sequences of perverse sheaves on $X$,
where the vertical arrows are induced by the
trace morphisms.
Hence by the induction on
the length of ${\cal F}$,
the assertion follows.
\qed
}

\begin{cor}\label{coralt}
Let $h\colon W\to X$ be a proper morphism of schemes
of finite type over a field $k$.
Let ${\cal F}$ be a perverse sheaf
on $X$
and ${\cal F}'$ be a perverse sheaf
on $W$.
Assume that there exists an open
subset $U\subset X$ satisfying the following conditions:
The base change $h_U\colon U'\to U$
of $h$ is finite and faithfully flat.
The perverse sheaf ${\cal F}$
is isomorphic to $j_{!*}{\cal F}_U$ for 
the open immersion $j\colon U\to X$
and ${\cal F}_U=j^*{\cal F}$.
The perverse sheaf
$\, ^{\rm p}\!H^0 h_U^*{\cal F}_U$ on $U'$
is isomorphic to a subquotient of the restriction
${\cal F}'|_{U'}$.

Then the perverse sheaf
${\cal F}$ on $X$
is isomorphic to a subquotient of 
$\, ^{\rm p}\!H^0 Rh_*{\cal F}'$.
\end{cor}

\proof{
Since $h_{U*}$ is $t$-exact,
we have
$h_{U*}({\cal F}'|_{U'})
=({}^{\rm p}\!H^0 Rh_*{\cal F}')|_U
$. Hence
$j_{!*}h_{U*}({\cal F}'|_{U'})
=j_{!*}({}^{\rm p}\!H^0 Rh_*{\cal F}')|_U$
is isomorphic to a subquotient of 
$\, ^{\rm p}\!H^0 Rh_*{\cal F}'$.
Since
$\, ^{\rm p}\!H^0 h_U^*{\cal F}_U$ on $U'$
is isomorphic to a subquotient of the restriction
${\cal F}'|_{U'}$,
the perverse sheaf
${\cal F}_U$ is a subquotient of 
$h_{U*}{\cal F}'|_{U'}$
by Proposition \ref{prffp}.
Hence ${\cal F}=j_{!*}{\cal F}_U$
is isomorphic to a subquotient of 
$j_{!*}(h_{U*}{\cal F}'|_{U'})$
and the assertion follows.
\qed
}

\medskip
Next, we consider the case
where $X$ is a scheme of
finite type over 
the spectrum $S$ of
a discrete valuation ring as in \cite[4.6]{au}.
Let $s$ and $\eta$ denote the
closed point and the generic point
of $S$ respectively
and let $i\colon X_s\to X$
and $j\colon X_\eta\to X$
be the closed immersion and the open immersion
of the fibers.
Let $\Lambda$ be
a finite field of characteristic $\ell$
invertible on $S$.
Then, we consider
the $t$-structure on
$D^b_c(X,\Lambda)$
obtained by gluing
(\cite[1.4.10]{BBD})
the $t$-structure
$(^{\rm p}\!D^{\leqq 0},\,
^{\rm p}\!D^{\geqq 0})$ on
$D^b_c(X_s,\Lambda)$ and
the $t$-structure
$(^{\rm p}\!D^{\leqq -1},\,
^{\rm p}\!D^{\geqq -1})$ on
$D^b_c(X_\eta,\Lambda)$.
In particular,
a constructible complex
${\cal F}\in D^b_c(X,\Lambda)$
is contained in $^{\rm p}\!D^{\leqq 0}$
if and only if we have
$i^*{\cal F}
\in\, ^{\rm p}\!D^{\leqq 0}$ and
$j^*{\cal F}\in\, ^{\rm p}\!D^{\leqq -1}$.

Note that if the $t$-structure
on $D^b_c(X_\eta,\Lambda)$
where $X_\eta$ is regarded
as a scheme over $\eta$ is
$(^{\rm p}\!D^{\leqq 0},\,
^{\rm p}\!D^{\geqq 0})$,
then that 
on $D^b_c(X_\eta,\Lambda)$
where $X_\eta$ is regarded
as a scheme over $S$ is
$(^{\rm p}\!D^{\leqq -1},\,
^{\rm p}\!D^{\geqq -1})$.
To distinguish them,
we call the former
the $t$-structure on
$X_\eta$ over $\eta$
and the latter
the $t$-structure on
$X_\eta$ over $S$.
We use the same
terminology for perverse sheaves on
$X_\eta$.

By the same argument as \cite[Th\'eor\`eme 4.3.1]{BBD}, 
we see that 
every object of
${\rm Perv}(X,\Lambda)$
is of finite length.
Further, simple objects are of the following form:
Let $V\subset X_s$ 
(resp.\ $V\subset X_\eta$) be an irreducible 
locally closed subset of dimension $d$
(resp.\ $d-1$) such that
the reduced part of the geometric fiber
$V_{\bar s}$ 
(resp.\ $V_{\bar \eta}$) is smooth.
Let $j\colon V\to X$ be the immersion
and let ${\cal G}$ be an irreducible
locally constant sheaf on $V$.
Then, $j_{!*}{\cal G}[d]$ is a simple perverse sheaf
on $X$.

The functors $j_!,Rj_*\colon
D^b_c(X_\eta,\Lambda)
\to D^b_c(X,\Lambda)$
are $t$-exact with respect to
the $t$-structure on
$X_\eta$ over $S$.
This follows from \cite[Th\'eor\`eme 3.1]{XIV}
by the argument in
\cite[4.6 (a)]{au}.
Let ${\cal F}
\in {\rm Perv}(X_\eta,\Lambda)$
be a perverse sheaf on
$X_\eta$ over $S$.
Then the intermediate extension
$j_{!*}{\cal F}\in 
{\rm Perv}(X,\Lambda)$
is defined as the image
$$j_{!*}{\cal F}
={\rm Im}(j_!{\cal F}
\to Rj_*{\cal F}).$$
Similarly as
\cite[(4.1.11.1)]{BBD},
the distinguished triangle
$j_!{\cal F}
\to
Rj_*{\cal F}
\to
i_*i^*Rj_*{\cal F}
\to$ 
and 
the $t$-exactness of
the functors $j_!$ and $Rj_*$
\cite[Appendix Remark (i)]{BBD},
imply the vanishing
${}^{\rm p}\!{\cal H}^qi^*
Rj_*{\cal F}=0$
for $q\neq 0,-1$ and define
an exact sequence
$$
0\to
i_*{}^{\rm p}\!{\cal H}^{-1}
i^*Rj_*{\cal F}
\to
j_!{\cal F}
\to
Rj_*{\cal F}
\to
i_*{}^{\rm p}\!{\cal H}^0
i^*Rj_*{\cal F}
\to 0$$
of perverse sheaves on $X$.
This induces
an isomorphism
\begin{equation}
{}^{\rm p}\!{\cal H}^{-1}
i^*Rj_*{\cal F}
\to
i^*{\rm Ker}(
j_!{\cal F}
\to
Rj_*{\cal F})
=
i^*{\rm Ker}(
j_!{\cal F}
\to
j_{!*}{\cal F})
\gets
i^*j_{!*}{\cal F}[-1]
\label{eqH0}
\end{equation}
of perverse sheaves on $X_s$
similarly as
\cite[(4.1.12.1)]{BBD}.

\subsection{Nearby cycles and local acyclicity}\label{ssS}

Assume that $S={\rm Spec}\ {\cal O}_K$ is the spectrum
of a strictly local discrete valuation ring.
Let $\bar \eta$ be a geometric point
above $\eta$.
Let $X$ be a scheme of finite type over $S$,
and let $i\colon X_s\to X$
and $j\colon X_\eta\to X$
denote the immersions
and let
$\pi\colon X_{\bar \eta}
\to X_\eta$
denote the canonical morphism.
Then, the nearby cycles
functor 
$$R\Psi= i^*R(j\pi)_*\pi^*
\colon D^b_c(X_\eta,\Lambda)
\to D^b_c(X_s,\Lambda)$$
is $t$-exact
with respect to
the $t$-structure
on $X_\eta$ over $\eta$
\cite[Corollaire 4.5]{au}.
Let $I={\rm Gal}(K_s/K)$
be the inertia group
and let $R\Gamma(I,-)$ denote
the derived functor for the
$I$-fixed part.
We have a canonical isomorphism
${\cal F}\to R\Gamma(I,\pi_*\pi^*{\cal F})$
and we identify
$i^*Rj_*{\cal F}=
R\Gamma(I,R\Psi {\cal F})$
by the induced isomorphism.

\begin{lm}\label{lmS}
Let $S={\rm Spec}\ {\cal O}_K$ be the spectrum
of a strictly local discrete valuation ring
and let $s$ and $\eta$ denote
the closed and the generic point of $S$
respectively.
Let $\bar \eta$ be the spectrum
of a separable closure $K_s$
of $K$
and $I={\rm Gal}(K_s/K)$
be the inertia group.
Let $X$ be a scheme of finite type over $S$,
and let $i\colon X_s\to X$
and $j\colon X_\eta\to X$
denote the immersions.
Let ${\cal F}$ be a perverse sheaf
of $\Lambda$-modules on $X_\eta$
{\em over} $S$.
Then the morphism
$i^*Rj_*{\cal F}\to R\Psi {\cal F}$
induces an isomorphism
\begin{equation}
i^*j_{!*}{\cal F}[-1]
\to (R\Psi {\cal F}[-1])^I
\label{eqI}
\end{equation}
to the inertia fixed part
as a perverse sheaf on
$X_s$.
\end{lm}

\proof{
Note that ${\cal F}[-1]$ is a perverse sheaf
on $X_\eta$ {\em over} $\eta$
and hence 
$R\Psi {\cal F}[-1]$ is a perverse sheaf on $X_s$
by the $t$-exactness of
$R\Psi$.
Let $P\subset I$ denote the
wild inertia subgroup.
Then, since the functor
taking the $P$-invariant parts
is an exact functor,
we have an isomorphism
$i^*Rj_*{\cal F}[-1]\to
R\Gamma(I/P,(R\Psi {\cal F}[-1])^P)$.
Since the profinite group
$I/P$ is cyclic, 
we have an isomorphism
$[\sigma-1\colon
(R\Psi {\cal F}[-1])^P\to (R\Psi {\cal F}[-1])^P]
\to
R\Gamma(I/P,(R\Psi {\cal F}[-1])^P)$
for a topological generator $\sigma$
of $I/P$.
Thus we obtain an isomorphism
$${}^{\rm p}\!{\cal H}^{-1}i^*
Rj_*{\cal F}\to (R\Psi {\cal F}[-1])^I$$
of perverse sheaves on $X_s$.
Hence the assertion follows from 
the isomorphism (\ref{eqH0}).
\qed}
\medskip

We study the local acyclicity
of a morphism to
the spectrum of a discrete valuation
ring with respect to a perverse sheaf.

\begin{pr}\label{prS}
Let $S={\rm Spec}\ {\cal O}_K$ be the spectrum
of a discrete valuation ring
and let $s$ and $\eta$ denote
the closed and the generic point of $S$
respectively.
Let $X$ be a scheme of finite type over $S$,
and let $i\colon X_s\to X$
and $j\colon X_\eta\to X$
denote the immersions.

{\rm 1.}
Let ${\cal G}$ be
a perverse sheaf of $\Lambda$-modules on $X$.
Assume that $X\to S$ is locally acyclic relatively
to ${\cal G}$.
Then ${\cal G}$ has no non-zero
subquotient supported
on the closed fiber
and is isomorphic to $j_{!*}j^*{\cal G}$.

{\rm 2.}
For a perverse sheaf ${\cal F}$ 
of $\Lambda$-modules on $X_\eta$
over $S$,
the following conditions are
equivalent:

{\rm (1)}
The morphism
$X\to S$ is locally acyclic relatively to
$j_{!*}{\cal F}$.

{\rm (2)}
Let $\bar s$ be a geometric point
above the closed point $s\in S$
and let $\bar i\colon X_{\bar s}\to X$
denote the canonical morphism.
Then, the canonical morphism
\begin{equation}
\bar i^*j_{!*}{\cal F}\to R\Psi {\cal F}
\label{eqij}
\end{equation}
is an isomorphism.

{\rm (3)}
The inertia group $I$ of $K$
acts trivially
on the nearby cycles complex $R\Psi{\cal F}$.

{\rm (4)}
The formation of
$j_{!*}{\cal F}$ commutes with the pull-back by
faithfully flat morphisms
$S'\to S$ of the spectra
of discrete valuation rings.
\end{pr}

\proof{
1.
We first show that
${\cal G}$ has no non-zero
subquotient supported
on the closed fiber.
The local acyclicity
is equivalent to
the vanishing 
$R\Phi {\cal G}=0$.
Since the shifted vanishing cycles functor
$R\Phi[-1]\colon
D^b_c(X,\Lambda)\to
D^b_c(X_s,\Lambda)$
is $t$-exact
\cite[Corollaire 4.6]{au},
it is reduced to the case
where ${\cal G}$ is
a simple perverse sheaf
by the induction on length of ${\cal G}$.
If ${\cal G}$ is supported
on the closed fiber,
we have $R\Phi{\cal G}[-1]
={\cal G}$.
Hence
${\cal G}$ has no non-zero
subquotient supported
on the closed fiber.

Since $j_{!*}j^*{\cal G}$
is the unique perverse sheaf
extending $j^*{\cal G}$
without non-trivial sub or quotient
perverse sheaf supported on the closed
fiber
by \cite[Corollaire 1.4.25]{BBD}, 
${\cal G}$ is canonically isomorphic to
$j_{!*}j^*{\cal G}$.

2.
(1)$\Leftrightarrow$(2):
The condition (2) is equivalent to
that for every geometric point
$x$ of $X_s$,
the canonical morphism
$j_{!*}{\cal F}_x\to 
R\Gamma(X_{(x)}\times_{S_{(s)}}\bar \eta,{\cal F})$ is an isomorphism.

(2)$\Leftrightarrow$(3):
Clear from the isomorphism
(\ref{eqI}).

(2)$\Rightarrow$(4):
Since the formation of nearby cycles complex
$R\Psi{\cal F}$
commutes with base change
\cite[Proposition 3.7]{TF},
the isomorphism (\ref{eqij})
implies the condition (4).

(4)$\Rightarrow$(2):
There exists a finite extension $K'$ of $K$
such that
the inertia action $I'\subset I$
on $R\Psi{\cal F}$
is trivial, since $\Lambda$ is a finite field.
Let $j'\colon X_{K'}\to X_{S'}$
denote the base change 
of the open immersion $j$
by $S'={\rm Spec}\, {\cal O}_{K'}\to S$,
let $i'\colon X_{\bar s}\to X_{S'}$
denote the canonical morphism 
and let ${\cal F}'$ denote the
pull-back of ${\cal F}$ on
$X_{K'}$.
We factorize the morphism (\ref{eqij})
as the composition of
$\bar i^*j_{!*}{\cal F}\to
\bar i^{\prime*}j'_{!*}{\cal F}'\to R\Psi {\cal F}.$
By (3)$\Rightarrow$(2) already proven,
the second arrow
is an isomorphism.
The condition (4) implies
that the first arrow is an isomorphism.
Hence the composition (\ref{eqij})
is an isomorphism.
\qed}\medskip

Finally, we consider the
case where
$X$ is a scheme of
finite type over 
a regular
noetherian connected
scheme $S$ of dimension $1$.
Let $\Lambda$ be
a finite field of characteristic $\ell$
invertible on $S$.
Then 
the $t$-structure
$(^{\rm p}\!D^{\leqq 0},\,
^{\rm p}\!D^{\geqq 0})$ on
$D^b_c(X,\Lambda)$ is defined
as the intersection of
the inverse images of
the $t$-structures
$(^{\rm p}\!D^{\leqq 0},\,
^{\rm p}\!D^{\geqq 0})$ on
$D^b_c(X\times_SS_s,\Lambda)$
for the base changes by
the localizations $S_s\to S$
at closed points $s\in S$.
If $Y=S$ is a smooth curve over a field
$k$ and if $f\colon X\to Y$
is a morphism of schemes
of finite type over $k$,
the $t$-structure on
$D^b_c(X,\Lambda)$ defined
above is the same as
that defined by considering
$X$ as a scheme of finite type over $k$.

\begin{cor}\label{corS}
Let $S$ be a regular
noetherian scheme of dimension $1$.
Let $X$ be a scheme of finite
type over $S$
and ${\cal F}$ be
a perverse sheaf of
$\Lambda$-modules on $X$.
Let $V\subset S$ be a dense
open subscheme such
that the base change
$X_V\to V$ is universally locally acyclic relatively to the
restriction ${\cal F}_V$ of ${\cal F}$.

Then, there exists a 
finite faithfully flat and generically \'etale morphism
$S'\to S$ of regular schemes
such that the base change
$X'\to S'$ is locally acyclic relatively to $j'_{!*}{\cal F}_{V'}$
where ${\cal F}_{V'}$
denotes the pull-back of ${\cal F}$
on $V'=V\times_SS'$
and $j'\colon X'_{V'}\to X'$
denotes the base change.
\end{cor}

\proof{
By Proposition \ref{prS}.2
(1)$\Rightarrow$(4)
and Lemma \ref{lmwa} below,
it suffices to consider locally
on a neighborhood of each point
of the complement $S\sm V$.
Since the coefficient field $\Lambda$
is finite, the assertion follows from
Proposition \ref{prS}.2 (3)$\Rightarrow$(1).
\qed}

\begin{lm}\label{lmwa}
Let $S$ be a regular
noetherian scheme of dimension $1$
and let $s_1,\ldots,s_n$ be closed points of $Y$.
Let $L_1,\ldots,L_n$ be finite separable extensions
of the local fields $K_1,\ldots,K_n$ of $S$ at $s_1,\ldots,s_n$.
Then, there exists a finite, faithfully flat and generically \'etale
morphism $S'\to S$
such that $S'\times_SK_i$ is isomorphic to
the disjoint union of finitely many copies of ${\rm Spec}\, L_i$.
\end{lm}

\proof{
Let $m$ be a common multiple of 
the degrees $[L_i:K_i]$
and let $A_i$ be the product of 
copies of $L_i$ such that $\dim_{K_i}A_i=m$.
Then, by weak approximation,
there exists a finite \'etale algebra $A$
over the fraction field $K$ of $S$
such that $A\otimes_KK_i=A_i$.
Then, it suffices to take the normalization $S'$ of
$S$ in $A$.
}

\subsection{$C$-transversality}
\label{ssCtr}

We introduce some terminology on 
proper intersection.

\begin{lm}\label{lmpi}
Let $f\colon C\to X$
and $h\colon W\to X$
be morphisms of schemes
of finite type over a field $k$.
Assume that $C$ is irreducible
of dimension $n$ and
that $h$ is locally of complete
intersection of relative 
virtual dimension $d$.
Then every irreducible component
of $h^*C=C\times_XW$
is of dimension $\geqq n+d$.
\end{lm}

\proof{
Since the assertion is local on $W$,
we may decompose 
$h=gi$ as the composition of
a smooth morphism $g$ with
a regular immersion of codimension
$c$.
Since the assertion is clear for
$g$, we may assume that
$h=i$ is a regular immersion.
Then, it follows from
\cite[Proposition (5.1.7)]{EGA4}.
\qed}

\begin{df}\label{dfpi}
Let $f\colon C\to X$
and $h\colon W\to X$
be morphisms of schemes
of finite type over a field $k$.
Assume that every irreducible component
of $C$ is of dimension $n$ and
that $h$ is locally of complete
intersection of relative 
virtual dimension $d$.
We say that $h\colon W\to X$
meets $f\colon C\to X$ {\em properly}
if $h^*C=C\times_XW$
is of dimension $n+d$.
\end{df}

By Lemma \ref{lmpi},
the condition that
$h^*C=C\times_XW$
is of dimension $n+d$
is equivalent to the condition that
every irreducible component
of $h^*C=C\times_XW$
is of dimension $n+d$.

\begin{lm}\label{lmpi2}
Let $f\colon C\to X$
be a morphism of schemes
of finite type over a field $k$.
Assume that $X$ is equidimensional of dimension $m$
and that
$C$ is equidimensional of dimension $n\geqq m$.
We consider the following conditions:

{\rm (1)}
Every morphism $h\colon W\to X$
locally of complete intersection
meets $C$ properly.

{\rm (2)}
For every closed point $x$ of $X$,
the fiber $C\times_Xx$
is of the dimension $n-m$.

{\rm 1.}
We have {\rm (2)}$\Rightarrow${\rm (1)}.
Assume that the condition {\rm (2)}
is satisfied
and let $h\colon W\to X$
be a morphism 
locally of complete intersection
of relative virtual dimension $d$
of schemes of 
finite type over $k$.
Then $C\times_XW$ is 
equidimensional of dimension $n+d$
and the morphism
$C\times_XW\to W$
satisfies the condition {\rm (2)}
and hence {\rm (1)}.

{\rm 2.}
If $X$ is regular,
we have
{\rm (1)}$\Rightarrow${\rm (2)}.

{\rm 3.}
Assume that $X={\mathbf P}$ is a projective space
and let $c$ be an integer.
Then, the linear subspaces
$V\subset {\mathbf P}$ of
codimension $c$
such that the immersion
$V\to {\mathbf P}$ meets $C$
properly form 
a dense open subset
of the Grassmannian variety
${\mathbf G}$.
\end{lm}

\proof{
1.
Assume that the condition (2)
is satisfied and 
let $h\colon W\to X$
be a morphism
locally of complete intersection
of relative dimension $d$.
Then, we have
$\dim C\times_XW\leqq \dim W+n-m
=n+d$.
Hence, $C\times_XW$
is equidimensional of dimension
$n+d$ by Lemma \ref{lmpi}.
The rest is clear.

2.
If $X$ is regular
and $x$ is a closed point,
the closed immersion
$i\colon x\to X$ is
a regular immersion of
codimension $m$ and hence
the condition (1) implies that
$\dim C\times_Xx=n-m$.

3.
Let ${\mathbf V}\subset
{\mathbf P}\times {\mathbf G}$
be the universal family
of linear subspaces of codimension
$c$
and we consider the cartesian diagram
$$\xymatrix{
&C_{\mathbf V}\ar[r]\ar[d]^{
\hspace{5mm}\square}\ar[ld]
&C\ar[d]\\
{\mathbf G}
&{\mathbf V}\ar[l]\ar[r]
&{\mathbf P}.
}$$
Then, since
the projection
${\mathbf V}\to {\mathbf P}$
is smooth of relative dimension
$\dim {\mathbf G}-c$,
we have
$\dim C_{\mathbf V}=
\dim {\mathbf G}+n-c$.
Hence 
the open subset
of ${\mathbf G}$
consisting of $V$
such that
$\dim C\times_{\mathbf P}V\leqq n-c$
is dense.
\qed}
\medskip

Recall that a closed subset
$C$ of a vector bundle
$E$ on a scheme $X$ is
said to be {\em conical} if
it is stable under the action
of the multiplicative group. 
For a closed conical subset
$C\subset E$,
the intersection $B=C\cap X$
with the $0$-section
identified with a closed subset
of $X$ is called the base of $C$.
We say that a morphism
$f\colon X\to Y$ of
noetherian schemes is finite
(resp. proper)
on a closed subset $Z\subset X$
if its restriction $Z\to Y$
is finite (resp. proper)
with respect to a
closed subscheme structure of $Z\subset X$.

\begin{df}\label{dfpt}
Let $f\colon X\to Y$ be
a morphism
of smooth schemes
over a field $k$ and
let $C\subset T^*X$ be a 
closed conical subset.

{\rm 1.
(\cite[1.2]{Be})}
We say that $f\colon X\to Y$
is $C$-{\em transversal}
if the inverse image of
$C$ by the canonical morphism
$X\times_YT^*Y\to T^*X$
is a subset of the $0$-section
$X\times_YT^*_YY
\subset X\times_YT^*Y$.

{\rm 2.}
Assume that every irreducible component
of $X$ is of dimension $n$
and that
every irreducible component
of $C$ is of dimension $n$.
Assume that every irreducible component
of $Y$ is of dimension $m\leqq n$.
We say that $f\colon X\to Y$
is {\em properly} $C$-transversal
if $f\colon X\to Y$
is $C$-transversal and if
for every closed point $y$ of
$Y$,
the fiber
$C\times_Yy$
is of dimension $n-m$.
\end{df}

\begin{df}\label{dfpth}
Let $h\colon W\to X$ be
a morphism
of smooth schemes
over a field $k$ and
let $C\subset T^*X$ be a 
closed conical subset.
Let $K\subset
W\times_XT^*X$ be the inverse image
of the $0$-section
$T^*_WW\subset T^*W$
by the canonical morphism
$W\times_XT^*X\to T^*W$.

{\rm 1.
(\cite[1.2]{Be})}
We say that $h\colon W\to X$
is $C$-{\em transversal}
if the intersection
$(W\times_XC)\cap K
\subset W\times_XT^*X$
is a subset of the $0$-section
$W\times_XT^*_XX$.

If $h\colon W\to X$ is
$C$-transversal,
a conical subset
$h^\circ C\subset T^*W$
is defined to be the image
of $h^*C=W\times_XC$
by $W\times_XT^*X\to T^*W$.

{\rm 2.
(\cite[Definition 7.1]{CC})}
Assume that every irreducible component
of $X$ is of dimension $n$
and that
every irreducible component
of $C$ is of dimension $n$.
Assume that every irreducible component
of $W$ is of dimension $m$.
We say that $h\colon W\to X$
is {\em properly} $C$-transversal
if $h\colon W\to X$
is $C$-transversal and if
$h\colon W\to X$
meets $C\to X$ properly.
\end{df}

If $h\colon W\to X$ is $C$-transversal,
the morphism
$W\times_XT^*X\to T^*W$
is finite on
$h^*C=W\times_XC$
and hence
$h^\circ C\subset T^*W$
is a {\em closed} subset 
by \cite[Lemma 1.2 (ii)]{Be}.
For a morphism
$r\colon X\to Y$
of smooth schemes
proper on the base
$B=C\cap T^*_XX\subset X$
of a closed conical subset
$C\subset T^*X$,
the closed conical subset
$r_\circ C\subset T^*Y$
is defined to be
the image by the projection
$X\times_YT^*Y\to T^*Y$
of the inverse image of $C$
by the canonical morphism
$X\times_YT^*Y\to T^*X$.

\begin{lm}\label{lmtrbc}
Let $f\colon X\to Y$ be
a {\em smooth} morphism
of smooth schemes
over a field $k$ and
let $C\subset T^*X$ be a 
closed conical subset.
Let
$$\begin{CD}
X@<h<< W\\
@VfVV
\hspace{-10mm}
\square
\hspace{7mm}
@VVgV\\
Y@<i<<Z
\end{CD}$$
be a cartesian diagram
of smooth schemes over
$k$.

{\rm 1.}
Assume that
$f\colon X\to Y$ is
$C$-transversal 
(resp.\ properly $C$-transversal).
Then,
$h\colon W\to X$ is
$C$-transversal 
(resp.\ properly $C$-transversal)
and
$g\colon W\to Z$ is
$h^\circ C$-transversal 
(resp.\ properly $h^\circ C$-transversal).

{\rm 2.}
Assume that
$f\colon X\to Y$ is proper
on the base of $C$.
Then, 
$i\colon Z\to Y$ is
$f_\circ C$-transversal
if and only if
$h\colon W\to X$ is
$C$-transversal.
If these equivalent conditions
are satisfied, we have
$i^\circ f_\circ C=
g_\circ h^\circ C$.
\end{lm}

\proof{
1.
The assertion for the transversality
is proved in \cite[Lemma 3.9.2]{CC}.
The proper transversality
of $h\colon W\to X$ follows from 
the transversality
and Lemma \ref{lmpi2} applied
to $C\to Y$ and $Z\to Y$.
The proper $h^\circ C$-transversality
of $g\colon W\to Z$ follows
from that
$h^*C\to h^\circ C$ is finite.

2.
We consider the commutative diagram
$$\begin{CD}
T^*X@<<< 
W\times_XT^*X@>{dh}>> T^*W\\
@AAA
\hspace{-15mm}
\square
\hspace{16mm}
@AAA
@AAA
\\
X\times_YT^*Y@<<< 
W\times_YT^*Y@>{g^*(di)}>> 
W\times_ZT^*Z\\
@VVV
\hspace{-15mm}
\square
\hspace{16mm}
@VV{g_*}V
\hspace{-17mm}
\square
\hspace{14mm}
@VVV
\\
T^*Y@<<< 
Z\times_YT^*Y@>{di}>> T^*Z
\end{CD}$$
with cartesian squares
indicated by $\square$.
The upper vertical arrows
are injections.
Since $dh$ induces an isomorphism
$W\times_XT^*X/Y\to T^*W/Z$
for the relative cotangent bundles
and $f\colon X\to Y$ is smooth,
the upper right square is
also cartesian.

Let $K$ and $K'$
be the inverse image
of the $0$-sections
by $dh\colon W\times_XT^*X\to T^*W$
and
$di\colon Z\times_YT^*Y\to T^*Z$
respectively.
Since the upper right square is cartesian,
$K$ is identified with
the inverse image of
the $0$-section by 
$g^*(di)\colon
W\times_YT^*Y
\to W\times_ZT^*Z$
which equals 
$g_*^{-1}(K')
\subset W\times_YT^*Y$.

Since the lower left square
is cartesian, 
the pull-back
$Z\times_Yf_\circ C$
is the image $g_*(C')$ of
$C'=(W\times_X C)
\cap 
(W\times_YT^*Y)$.
Hence the condition that
$(Z\times_Yf_\circ C) \cap K'
=g_*(C')\cap K'
=g_*(C'\cap g_*^{-1}(K'))$ 
is a subset
of the $0$-section
is equivalent to the condition
that $(W\times_X C)
\cap K=C'\cap g_*^{-1}(K')$ is a subset
of the $0$-section.

If these conditions are satisfied,
the equality
$i^\circ f_\circ C=
g_\circ h^\circ C$
follows from the cartesian diagram.
\qed}

\begin{lm}\label{lmLeg}
Let ${\mathbf P}$
be a projective space of
dimension $n$ and
let $C\subset T^*{\mathbf P}$
be a closed conical subset.

{\rm 1.}
Let ${\mathbf P}^\vee$
be the dual projective space,
let
$Q\subset {\mathbf P}
\times {\mathbf P}^\vee$
be the universal family of
hyperplanes
and let 
\begin{equation}
\begin{CD}
{\mathbf P}
@<p<< Q
@>{p^\vee}>>
{\mathbf P}^\vee
\end{CD}
\label{eqpv}
\end{equation}
be the projections.
Let $C^\vee
=p^\vee_\circ p^\circ C$
be the Legendre transform.
Let $V\subset {\mathbf P}$
be a linear subspace
and let $V^\vee\subset {\mathbf P}^\vee$
be the dual subspace.
Then the immersion
$V\to {\mathbf P}$
is $C$-transversal if and only if
$V^\vee\to {\mathbf P}^\vee$
is $C^\vee$-transversal.

{\rm 2.}
Assume that 
every irreducible component of
$C$ is of dimension $n
=\dim {\mathbf P}$ and
let $0\leqq c\leqq n$ be an integer.
Then, the linear subspaces
$V\subset {\mathbf P}$ of
codimension $c$
such that the immersion
$V\to {\mathbf P}$ is
properly $C$-transversal form 
a dense open subset
of the Grassmannian variety
${\mathbf G}$.
\end{lm}

\proof{
1.
The $C$-transversality of
$V\to {\mathbf P}$ means
${\mathbf P}(T^*_V{\mathbf P})
\cap {\mathbf P}(C)=\varnothing$
and similarly for
the $C^\vee$-transversality of
$V^\vee\to {\mathbf P}^\vee$.
Then, the assertion follows
from 
${\mathbf P}(T^*_V{\mathbf P})
={\mathbf P}(T^*_{V^\vee}{\mathbf P}^\vee)$
and ${\mathbf P}(C)=
{\mathbf P}(C^\vee)$
under the identification
${\mathbf P}(T^*{\mathbf P})
=Q
={\mathbf P}(T^*{\mathbf P}^\vee)$.

2.
First, we show that the condition on $V$ is an open condition.
For the transversality, it suffices to apply a relative version of
\cite[Lemma 3.4.4]{CC} or equivalently
\cite[Lemma 1.2(i)]{Be} to
the closed immersion of the universal
family ${\mathbf V}\to {\mathbf P}\times {\mathbf G}$
since  the projection ${\mathbf V}\to {\mathbf G}$
is proper.
For the proper transversality,
it follows from the semi-continuity of fiber dimension
\cite[Th\'eor\`eme (13.1.3)]{EGA4}

It suffices to show the existence of $V$.
By induction on $c$,
it is reduced to the case $c=1$.
By 1,
the hyperplanes $H$ such that
the immersion $H\to {\mathbf P}$
is $C$-transversal
is parametrized by
the complement of
the image $p^\vee({\mathbf P}(C))
\subsetneqq {\mathbf P}^\vee$.
Hence, the assertion follows
from this and Lemma \ref{lmpi2}.3.
\qed}

\subsection{$SS{\cal F}$-transversality}
\label{ssFtr}

For the definitions and
basic properties 
of the singular support
of a constructible complex
on a smooth scheme over a field,
we refer to \cite{Be}
and \cite{CC}.
Let $k$ be a field and 
let $\Lambda$ be a finite field
of characteristic $\ell$ invertible in $k$.
Let $X$ be a smooth scheme over $k$
such that every irreducible
component is of dimension $n$
and let
${\cal F}$ be a constructible complex on $X$.
The singular support
$SS{\cal F}$ is defined in \cite{Be}
as a closed conical subset of
the cotangent bundle $T^*X$.
By \cite[Theorem 1.3 (ii)]{Be},
every irreducible component $C_a$
of the singular support
$$SS{\cal F}=C=\bigcup_aC_a$$
is of dimension $n=\dim X$.

Further if $k$ is perfect,
the characteristic cycle
$$CC{\cal F}=\sum_am_aC_a$$
is defined as a linear combination
with ${\mathbf Z}$-coefficients
in \cite[Definition 5.10]{CC}.
It is characterized by the
Milnor formula
\begin{equation}
-\dim {\rm tot}
\phi_u({\cal F},f)=
(CC{\cal F},df)_{T^*U,u}
\label{eqMil}
\end{equation}
for morphisms $f\colon U\to Y$
to smooth curves $Y$
defined on an \'etale neighborhood $U$
of an isolated characteristic point $u$.
For more detail on the notation, 
we refer to \cite[Section 5.2]{CC}.

\begin{lm}\label{lmh}
Let $h\colon W\to X$
be a morphism of smooth
schemes over a field $k$.
Let ${\cal F}$
be a constructible complex
of $\Lambda$-modules on $X$
and let $C$ denote the singular
support $SS{\cal F}$.
If $h\colon W\to X$
is properly $C$-transversal,
we have 
$$SSh^*{\cal F}=
h^\circ SS{\cal F}.$$
\end{lm}

\proof{
By \cite[Theorem 1.4 (iii)]{Be},
we may assume that $k$
is perfect.
Suppose $\dim W=\dim X+d$. 
If ${\cal F}$ is a perverse sheaf on $X$,
then $h^*{\cal F}[d]$ is a perverse sheaf on $W$
by the assumption that
$h\colon W\to X$
is $C$-transversal and by
\cite[Lemma 8.6.5]{CC}.
Hence by \cite[Theorem 1.4 (ii)]{Be},
we may assume that 
${\cal F}$ is a perverse sheaf.
By \cite[Proposition 5.14.2]{CC},
we have $CC{\cal F}\geqq 0$ and
the support of $CC{\cal F}$
equals the singular support
$SS{\cal F}$.
Also
we have $(-1)^dCCh^*{\cal F}\geqq 0$ and
the support of $CCh^*{\cal F}$
equals the singular support
$SSh^*{\cal F}$.

By the assumption that
$h\colon W\to X$
is properly $C$-transversal and
by \cite[Theorem 7.6]{CC},
we have $CCh^*{\cal F}
=h^!CC{\cal F}$.
Hence by the positivity
\cite[Proposition 7.1 (a)]{Ful}, the singular support
$SSh^*{\cal F}$
equals the support $h^\circ SS{\cal F}$
of $h^!CC{\cal F}$.
\qed}

\begin{lm}\label{lmtrans}
Let $k$ be a field and
$f\colon X\to Y$ be
a morphism of schemes 
of finite type over $k$.
Assume that
$Y$ is smooth over $k$.
Let ${\cal F}$
be a constructible complex 
of $\Lambda$-modules on $X$.
Let 
$$\begin{CD}
X@>i>> P\\
@V{i'}VV @VVgV\\
P'@>{g'}>>Y
\end{CD}$$
be a commutative
diagram of schemes over $k$
such that
$i$ and $i'$ are closed immersions and
the schemes $P$ and $P'$
are smooth over $k$.
Let
$C=SSi_*{\cal F}\subset
X\times_PT^*P\subset T^*P$
and
$C'=SSi'_*{\cal F}\subset
X\times_{P'}T^*P'\subset T^*P'$
denote the singular supports
of the direct images.
Then,
$P\to Y$ is $C$-transversal 
(resp.\ properly $C$-transversal)
if and only if
$P'\to Y$ is $C'$-transversal 
(resp.\ properly $C'$-transversal).
\end{lm}

\proof{
By factorizing $P\to Y$
as the composition of
the graph $P\to P\times Y$
and the projection
$P\times Y$,
we may assume that
$P\to Y$ is smooth.
Similarly,
we may assume that
$P'\to Y$ is smooth.
By considering the morphism
$(i, i')\colon X\to
P\times_YP'$,
we may assume that
there exists a smooth
morphism
$P'\to P$ compatible
with the immersions of $X$
and the morphisms to $Y$.
Since the assertion is \'etale local
on $P$, we may assume
that there exists
a section $s\colon P\to P'$
compatible with the immersions of $X$
and the morphisms to $Y$.
Then, we have
$C'=s_\circ C$
and the assertion follows
from \cite[Lemma 3.8]{CC}.
\qed}
\medskip

Lemma \ref{lmtrans}
allows us to make the
following definition.

\begin{df}\label{dftrans}
Let $k$ be a field and
$f\colon X\to Y$ be
a morphism of schemes 
of finite type over $k$.
Assume that
$Y$ is smooth over $k$
but we do {\em not} require $X$ to be smooth.
Let ${\cal F}$
be a constructible complex 
of $\Lambda$-modules on $X$.

We say that $f\colon X\to Y$
is $SS{\cal F}$-transversal
(resp.\ properly 
$SS{\cal F}$-transversal)
if locally on $X$
there exist a closed immersion
$i\colon X\to P$
to a smooth scheme $P$ over $k$
and a morphism
$g\colon P\to Y$ over
$k$ such that $f=g\circ i$
and that
$g\colon P\to Y$ 
is $C$-transversal
(resp.\ properly 
$C$-transversal)
for $C=SSi_*{\cal F}$.
\end{df}

In Definition \ref{dftrans},
we obtain an equivalent
condition by requiring
that $g$ is {\em smooth}.

Let $f\colon X\to Y$ be
a morphism of schemes 
of finite type over a field $k$
such that
$Y$ is smooth over $k$
and let ${\cal F}$
be a constructible complex 
of $\Lambda$-modules on $X$.
For an open subset
$U\subset X$,
we say $f\colon X\to Y$
is $SS{\cal F}$-transversal
(resp.\ properly 
$SS{\cal F}$-transversal)
on $U$ 
if the restriction 
$U\to Y$ of $f$
is $SS{\cal F}_U$-transversal
(resp.\ properly 
$SS{\cal F}_U$-transversal)
for the restriction
${\cal F}_U$ of ${\cal F}$ on $U$.
Similarly, for an open subset
$V\subset Y$,
we say $f\colon X\to Y$
is $SS{\cal F}$-transversal
(resp.\ properly 
$SS{\cal F}$-transversal)
on $V$ 
if the base change
$X\times_YV\to V$ of $f$
is $SS{\cal F}_{X\times_YV}$-transversal
(resp.\ properly 
$SS{\cal F}_{X\times_YV}$-transversal)
for the restriction
${\cal F}_{X\times_YV}$ of ${\cal F}$ on $X\times_YV$.

\begin{lm}\label{lmtr}
Let $f\colon X\to Y$ be
a morphism of schemes 
of finite type over a field $k$.
Assume that
$Y$ is smooth over $k$.
Let ${\cal F}$
be a constructible complex 
of $\Lambda$-modules on $X$.

{\rm 1.}
Assume that $f\colon X\to Y$ is smooth
and that ${\cal F}$ is locally constant.
Then, $f\colon X\to Y$ is
properly $SS{\cal F}$-transversal.

{\rm 2.}
Assume that 
$f\colon X\to Y$ is
$SS{\cal F}$-transversal.
Or more weakly,
suppose that
there exists a quasi-finite faithfully
flat morphism 
$Y'\to Y$ of smooth schemes
over $k$ such that
the base change
$f'\colon X'\to Y'$ is $SS{\cal F}'$-transversal
for the pull-back ${\cal F}'$
of ${\cal F}$ on $X'=X\times_YY'$.

Then,
$f\colon X\to Y$ is
universally locally acyclic
relatively to 
${\cal F}$.

{\rm 3.}
The following conditions are equivalent:

{\rm (1)}
$f\colon X\to Y$ is
$SS{\cal F}$-transversal
(resp.\ properly 
$SS{\cal F}$-transversal).

{\rm (2)}
For every integer $q$
and for every constituant ${\cal G}$
of the perverse sheaf $^{\rm p}\!{\cal H}^q{\cal F}$,
the morphism
$f\colon X\to Y$ is
$SS{\cal G}$-transversal
(resp.\ properly 
$SS{\cal G}$-transversal).

{\rm 4.}
Let $k'$ be an extension of $k$.
Then $f\colon X\to Y$ is
$SS{\cal F}$-transversal
(resp.\ properly 
$SS{\cal F}$-transversal)
if and only if
the base change $f'\colon X'\to Y'$ 
by ${\rm Spec}\ k'\to {\rm Spec}\ k$ is
$SS{\cal F}'$-transversal
(resp.\ properly 
$SS{\cal F}'$-transversal)
for the pull-back ${\cal F}'$
on $X'$ of ${\cal F}$.
\end{lm}

\proof{
1.
If ${\cal F}$ is locally constant,
then the singular support
$SS{\cal F}$ is a subset of the $0$-section
$T^*_XX$.
Hence the assertion follows.

Since the remaining
assertions 2-4 are local on $X$,
we may take a closed
immersion $i\colon X\to P$
to a smooth scheme $P$ over $k$
such that $f$ is the composition
of $i$ with a morphism
$P\to Y$ over $k$.
Replacing $X$ and ${\cal F}$ 
by $P$ and $i_*{\cal F}$,
we may assume that $X$
is smooth over $k$.
Set $C=SS{\cal F}$.

2.
If $f\colon X\to Y$
is $C$-transversal,
the morphism
$f\colon X\to Y$ is 
universally locally acyclic
relatively to 
${\cal F}$
by the definition of
singular support.
Thus under the weaker assumption,
the morphism
$f'\colon X'\to Y'$ is universally 
locally acyclic with respect
to the pull-back ${\cal F}'$.
Since $Y'\to Y$ is quasi-finite
and faithfully flat
and since the local acyclicity
descends for faithfully flat morphisms,
the morphism
$f\colon X\to Y$ itself is universally 
locally acyclic with respect
to ${\cal F}$.

3.
By \cite[Theorem 1.4 (ii)]{Be},
the singular support
$SS{\cal F}$ equals the union
of $SS{\cal G}$
for the constituants ${\cal G}$
of the perverse sheaves 
$^{\rm p}\!{\cal H}^q{\cal F}$
for integers $q$.
Hence the assertion follows.

4.
By \cite[Theorem 1.4 (iii)]{Be},
the construction of the singular support
commutes with change of
base fields.
Hence the assertion follows.
\qed}

\begin{lm}\label{lmtrV}
Let $f\colon X\to Y$ be
a morphism of schemes 
of finite type over a field $k$.
Assume that
$Y$ is smooth over $k$.
Let ${\cal F}$
be a constructible complex 
of $\Lambda$-modules on $X$.
Assume that $f\colon X\to Y$ is
$SS{\cal F}$-transversal.

{\rm 1.}
Assume that ${\cal F}$
is a perverse sheaf.
Let $V\subset Y$
be a dense open subscheme
and 
$j\colon X_V=X\times_YV\to X$
be the open immersion.
Then, there is a unique isomorphism
${\cal F}\to j_{!*}j^*{\cal F}$
such that the restriction on $X_V$ is
the identity.

{\rm 2.}
There exists
a dense open subscheme
$V\subset Y$ such that
the base change $f\colon X_V\to V$ is
properly $SS{\cal F}$-transversal
on $V$.
\end{lm}

\proof{
1.
By \cite[Corollaire 1.4.25]{BBD},
it suffices to show
that for every constituant of ${\cal F}$,
its restriction on $X_V$ is non-trivial.
Let ${\cal G}$ be a constituant of ${\cal F}$.
By Lemma \ref{lmtr}.3 and 2,
the morphism
$f\colon X\to Y$ is
locally acyclic relatively to ${\cal G}$.
Let $x$ be a geometric point of $X$ 
such that ${\cal G}_x\neq 0$
and let 
$y\to f(x)$ be a specialization
for a geometric point $y$ of $V$.
Then, since the canonical morphism
${\cal G}_x\to
R\Gamma(X_{(x)}\times_{Y_{(f(x))}}y,{\cal G})$
is an isomorphism,
the restriction of ${\cal G}$ on $X_V$
is non-trivial. Thus the assertion is proved.

2.
As in the proof of Lemma \ref{lmtr},
we may assume that $X$
is smooth over $k$.
Set $C=SS{\cal F}$.
There exists a dense open subset
$V\subset Y$
such that for every irreducible
component $C_a$ with 
the reduced scheme structure
of $C=\bigcup_aC_a$,
the base change
$C_a\times_YV\to V$ is flat.
\qed}

\begin{lm}\label{lmtrb}
Let $f\colon X\to Y$ be
a morphism of schemes 
of finite type over a field $k$.
Assume that
$Y$ is smooth over $k$.
Let ${\cal F}$
be a constructible complex 
of $\Lambda$-modules on $X$.
Let $Y'\to Y$ be a morphism
of smooth schemes over $k$
and let
$$\begin{CD}
X@<h<< X'\\
@VfVV 
\hspace{-10mm}
\square
\hspace{7mm}
@VV{f'}V\\
Y@<<< Y',
\end{CD}$$
be a cartesian diagram.
Let ${\cal F}'$ denote the
pull-back of ${\cal F}$ on $X'$.

{\rm 1.}
We consider the following conditions:

{\rm (1)}
$f\colon X\to Y$ is
$SS{\cal F}$-transversal
(resp.\ properly 
$SS{\cal F}$-transversal).

{\rm (2)}
$f'\colon X'\to Y'$ is
$SS{\cal F}'$-transversal
(resp.\ properly 
$SS{\cal F}'$-transversal).

Then, we have
{\rm (1)}$\Rightarrow${\rm (2)}.
Conversely,
if $Y'\to Y$ is \'etale surjective,
we have
{\rm (2)}$\Rightarrow${\rm (1)}.

{\rm 2.}
Assume that $f\colon X\to Y$ is
$SS{\cal F}$-transversal, that
${\cal F}$ is a perverse sheaf on $X$
and that $\dim Y'=\dim Y+d$.
Then ${\cal F}'[d]$ is
a perverse sheaf on $X'$.

{\rm 3.}
Assume that $f\colon X\to Y$ is smooth
and is properly $SS{\cal F}$-transversal.
Then, we have
$SS{\cal F}'=h^\circ SS{\cal F}$.
Further if $k$ is perfect,
we have
$CC{\cal F}'=h^! CC{\cal F}$.
\end{lm}

\proof{
Since the assertions are local on $X$,
we may take a closed
immersion $i\colon X\to P$
to a smooth scheme $P$ over $Y$.
As in the proof of Lemma \ref{lmtr},
we may assume that $f\colon
X\to Y$ is smooth.
Set $C=SS{\cal F}$.

1.
Assume that $f\colon X\to Y$
is $C$-transversal.
The pair $(h,f')$ of morphisms
is $C$-transversal by
Lemma \ref{lmtrbc}.1.
Hence, 
${\cal F}'=h^*{\cal F}$
is micro-supported on $h^\circ C$
by \cite[Lemma 4.2.4]{CC}
and we have
an inclusion
$SS{\cal F}'\subset h^\circ C$
and $f'$ is $SS{\cal F}'$-transversal.
Thus the implication 
(1)$\Rightarrow$(2) is proved
for the $C$-transversality.
The assertion on
the proper $C$-transversality
follows from this
and Lemma \ref{lmtrbc}.1.

Since the formation of
singular support is \'etale local,
we have
(2)$\Rightarrow$(1) 
if $Y'\to Y$ is \'etale surjective.

{\rm 2.}
Since $h\colon X'\to X$
is $C$-transversal by Lemma \ref{lmtrbc}.1,
the assertion follows from
\cite[Lemma 8.6.5]{CC}.

{\rm 3.}
Since $h\colon X'\to X$
is properly $C$-transversal by Lemma \ref{lmtrbc}.1,
the assertion for
$SS{\cal F}'$
(resp.\ for $CC{\cal F}'$)
follows from 
Lemma \ref{lmh}
(resp.\ \cite[Theorem 7.6]{CC}).
\qed}

\medskip
Lemma \ref{lmtrb}.3 is
closely related to the subject studied in
\cite{EH}.

\begin{lm}\label{lmtrfun}
Let $f\colon X\to Y$ be
a morphism of schemes 
of finite type over a field $k$.
Assume that
$Y$ is smooth over $k$.
Let ${\cal F}$
be a constructible complex 
of $\Lambda$-modules on $X$.

{\rm 1.}
Let $g\colon Y\to Z$
be a smooth morphism
of smooth schemes over $k$.
If $f\colon X\to Y$ is
$SS{\cal F}$-transversal
(resp.\ properly 
$SS{\cal F}$-transversal),
then the composition
$gf\colon X\to Z$ is
$SS{\cal F}$-transversal
(resp.\ properly 
$SS{\cal F}$-transversal).

{\rm 2.}
Let
$h\colon W\to X$ be a smooth
morphism
of schemes of finite type over $k$.
If $f\colon X\to Y$ is
$SS{\cal F}$-transversal
(resp.\ properly 
$SS{\cal F}$-transversal),
then the composition
$fh\colon W\to Y$ is
$SSh^*{\cal F}$-transversal
(resp.\ properly 
$SSh^*{\cal F}$-transversal).

{\rm 3.}
Let $$
\xymatrix{
X\ar[r]^r\ar[rd]_f&
X'\ar[d]^{f'}\\
&Y}$$
be a commutative diagram
of morphisms
of schemes of finite type over $k$.
Assume that
$r\colon X\to X'$ is {\rm proper}
on the support of ${\cal F}$
and that $f\colon X\to Y$ is 
quasi-projective.
If $f\colon X\to Y$ is
$SS{\cal F}$-transversal,
then $f'\colon X'\to Y$ is
$SS\, Rr_*{\cal F}$-transversal.
\end{lm}

\proof{
1.
As in the proof of Lemma \ref{lmtr},
we may assume that $X$
is smooth over $k$.
Set $C=SS{\cal F}$.
Since $g\colon Y\to Z$
is smooth,
the $C$-transversality of
$f$ implies that of $gf$
by \cite[Lemma 3.6.3]{CC}.
The assertion on
the proper $C$-transversality
follows from this and
the smoothness of $g$.

2. 
Since the question is \'etale local
on $W$,
we may assume that
there exists a cartesian diagram
$$\xymatrix{
W\ar[r]^h\ar[d]^{
\hspace{5mm}\square}&
X\ar[r]^f\ar[d]^i&Y\\
Q\ar[r]&P\ar[ru]
}$$
of morphisms of schemes over $k$
such that
the vertical arrows are
closed immersions 
and the horizontal arrow
$Q\to P$ 
is a smooth morphism of smooth schemes over $k$.
By replacing $X$ and ${\cal F}$
by $P$ and $i_*{\cal F}$,
we may assume that $X$ is smooth.
Since
$W\times_XT^*X\to T^*W$
is an injection and
$SSh^*{\cal F}
=h^\circ SS{\cal F}$
by Lemma \ref{lmh},
the assertion follows.

3.
Since the assertion is
local on $X'$,
we may assume that
$X'$ and $Y$ are affine
and hence $X$ is quasi-projective.
By taking
a closed immersion $i'\colon
X'\to P'$ to an affine space
and by factorizing
$X'\to Y$
as the composition of
the immersion
$(i',f')\colon X'\to P'\times Y$
and the projection
$P'\times Y\to Y$,
we may assume that
$X'$ is smooth.
Similarly, we take
an open subscheme $P$
of a projective
space and a closed immersion
$i\colon X\to P$.
Then, by factorizing
$X\to X'$
as the composition of
the immersion
$(i,r)\colon X\to P\times X'$
and the projection
$P\times X'\to X'$,
we may also assume that
$X$ is smooth,
by \cite[Lemma 3.8
(2)$\Rightarrow$(1)]{CC}.
By \cite[Lemma 2.2 (ii)]{Be},
we have $SS\, Rr_*{\cal F}
\subset r_\circ SS{\cal F}$.
Hence the assertion follows
from \cite[Lemma 3.8 (2)$\Rightarrow$(1)]{CC}.
\qed}

\medskip

We give two methods to
establish $SS{\cal F}$-transversality.

\begin{lm}\label{lmtrZ}
Let $Y\to S$ be
a smooth morphism 
of smooth schemes 
of finite type over a field $k$
and let $f\colon X\to Y$ be
a morphism of schemes 
of finite type over a field $k$.
Let ${\cal F}$
be a constructible complex 
of $\Lambda$-modules on $X$.
Assume that the composition
$X\to S$ is {\em properly}
$SS{\cal F}$-transversal.

{\rm 1.}
Assume that $k$ is perfect.
Then, the following conditions
are equivalent:

{\rm (1)}
$f\colon X\to Y$ is $SS{\cal F}$-transversal
(resp.\ properly
$SS{\cal F}$-transversal).

{\rm (2)}
For every closed point $s\in S$,
the fiber
$f_s\colon X_s\to Y_s$ is 
$SS{\cal F}_s$-transversal
(resp.\ properly
$SS{\cal F}_s$-transversal)
for the pull-back
${\cal F}_s$
of ${\cal F}$ on
$X_s=X\times_Ss$.

{\rm 2.}
Assume that
${\cal F}$ is a perverse sheaf on $X$
and that 
$f\colon X\to Y$ is
locally acyclic relatively to 
${\cal F}$.
If there exists a closed subset
$Z\subset X$ 
{\em quasi-finite} over $S$ such that
$f\colon X\to Y$
is $SS{\cal F}$-transversal
(resp.\ properly 
$SS{\cal F}$-transversal)
on the complement
of $Z$, then
$f\colon X\to Y$
is $SS{\cal F}$-transversal
(resp.\ properly 
$SS{\cal F}$-transversal)
on $X$.
\end{lm}

\proof{
1.
The implication (1)$\Rightarrow$(2)
is a special case of
Lemma \ref{lmtrb}.1.
We show
(2)$\Rightarrow$(1).
Since the question is local on $X$,
we may assume that $f\colon X\to Y$
is smooth.
Let $T^*X/S$ and
$T^*Y/S$ denote
the relative cotangent bundles
and let
$C=SS{\cal F}$.
By the assumption that
$X\to S$ is $C$-transversal,
the canonical surjection
$T^*X\to T^*X/S$
is finite on $C$
by \cite[Lemma 1.2 (ii)]{Be}.
Hence its image
$\bar C
\subset T^*X/S$
is a closed conical subset
and $C\to \bar C$ is finite.

The morphism
$X\to Y$ is $C$-transversal
if and only if
the inverse image of $\bar C
\subset T^*X/S$
by the canonical injection
$X\times_YT^*Y/S\to T^*X/S$
is a subset of the $0$-section.
This is equivalent to that
for every closed point 
$s\in S$ and
the closed immersion
$i_s\colon X_s\to X$,
the morphism
$f_s\colon X_s\to Y_s$
is $i_s^\circ C$-transversal.
Further, under the
assumption that
$f\colon X\to Y$ is $C$-transversal,
this is properly $C$-transversal
if and only if
$f_s\colon X_s\to Y_s$
is properly $i_s^\circ C$-transversal
for every closed point 
$s\in S$.

By the assumption that
$X\to S$ is properly $SS{\cal F}$-transversal and by
Lemma \ref{lmtrb}.3,
we have 
$SS{\cal F}_s=
i_s^\circ SS{\cal F}=
i_s^\circ C$
for every closed point
$s\in S$.
Hence the assertion is proved.

2.
By Lemma \ref{lmtr}.4,
we may assume that $k$
is perfect.
By 1 and Lemma \ref{lmtrb}.2, we may assume
that $S$ is a point
and further that
$S={\rm Spec}\ k$.
As in the proof of Lemma \ref{lmtrb},
after replacing $X$ by a smooth scheme $P$
over $Y$ containing $X$ as a closed subscheme,
we may assume that $f\colon X\to Y$
is smooth of relative dimension $d$.
Let $u\in Z$.
By replacing $X$ 
by a neighborhood of $u$,
we may assume $Z=\{u\}$.
Set $C=SS{\cal F},\ v=f(u)\in Y$
and regard
$X\times_YT^*Y$ as a closed subscheme
of $T^*X$.

We show that
$f\colon X\to Y$ is $C$-transversal,
assuming that $f\colon X\to Y$
is locally acyclic relatively to
${\cal F}$.
Namely, we show that
the intersection $C'
=C\cap (X\times_YT^*Y)$ is 
a subset of the $0$-section
$X\times_YT^*_YY$.
By the assumption that
$f\colon X\to Y$
is $C$-transversal outside $u$,
the intersection $C'
=C\cap (X\times_YT^*Y)$ is 
a subset of the union
$(X\times_YT^*_YY)
\cup (u\times_YT^*Y)$  with
the fiber at $u$.
Let $v=f(u)$ and
$\omega
\in u\times_YT^*Y
=v\times_YT^*Y$
be a non-zero element.
After shrinking $Y$
to a neighborhood of $v$
if necessary,
we take a smooth morphism
$Y\to {\mathbf A}^1
={\rm Spec}\ k[t]$
such that $dt(v)=\omega$.
Then, by
\cite[Lemma 3.6.3]{CC},
on a neighborhood of $u$,
the composition
$g\colon X\to Y\to {\mathbf A}^1$
is $C$-transversal except at most at $u$.
In other words, the point
$u$ is at most an isolated $C$-characteristic point
(\cite[Definition 5.3.1]{CC})
of 
$g\colon X\to {\mathbf A}^1$.

Since ${\cal F}$ is a perverse sheaf,
the characteristic cycle
$CC{\cal F}$ is an effective
cycle and its support equals
$C=SS{\cal F}$
by 
\cite[Proposition 5.14]{CC}.
Let $dg$ denote the section of
$X\times_YT^*Y\subset T^*X$ 
defined by the function
$g^*(t)$.
Since the composition 
$X\to Y\to {\mathbf A}^1$
is locally acyclic relatively to
${\cal F}$ by
\cite[Corollaire 5.2.7]{app},
we have
$(CC{\cal F},dg)_{T^*X,u}=0$
by the Milnor formula (\ref{eqMil}).
Therefore by the positivity
\cite[Proposition 7.1 (a)]{Ful}, the intersection
$SS{\cal F}\cap dg
=C'\cap dg$
is empty and
hence $\omega\notin C'$. 
Since $\omega$ is any non-zero
element of $u\times_YT^*Y$,
we conclude that
$C'\cap 
(u\times_YT^*Y)\subset 0$
and that
$f\colon X\to Y$
is $C$-transversal.

Assume further that
$f\colon X\to Y$
is properly $C$-transversal
outside $u$.
Since
$f\colon X\to Y$
is $C$-transversal,
the morphism
$T^*X\to T^*X/Y$
to the relative cotangent bundle
is finite on $C$
by \cite[Lemma 1.2 (ii)]{Be}
and the image
$\bar C\subset T^*X/Y$
of $C$ is a closed conical subset.
It is sufficient to
show that for every point
$y\in Y$,
the fiber $\bar C\times_Yy$
is of dimension $\leqq d$.
For $y\neq f(u)$,
this follows from the assumption. 
Assume $y=f(u)$.
Then, every irreducible component of
$\bar C\times_Yy$
is either 
a closure of an irreducible
component of
$\bar C\times_Yy
\cap (X\times_Yy\sm \{u\})$
or a subset of
the fiber
$T^*_u(X\times_Yy)$.
Since $\dim 
T^*_u(X\times_Yy)=d$,
the assertion is proved.
\qed}

\begin{lm}\label{lmtrqf}
Let 
$$\begin{CD}
W@>h>> X@>f>> Y\\
@A{j'}AA 
\hspace{-10mm}
\square
\hspace{7mm}
@AAjA\\
U'@>{h_U}>> U
\end{CD}$$
be a cartesian diagram of schemes 
of finite type over a field $k$.
Assume that $Y$ is smooth over $k$
and that $j$ is an open immersion.
Let ${\cal F}$ and ${\cal F}'$
be perverse sheaves
of $\Lambda$-modules on $X$
and on $W$ respectively
and let ${\cal F}_U$ and ${\cal F}'_{U'}$
be the restrictions on $U$
and on $U'$ respectively.
Assume that
${\cal F}$ is isomorphic to $j_{!*}{\cal F}_U$
and that
the perverse sheaf
${}^{\rm p}\!H^0(h_U^*{\cal F}_U)$
on $U'$ is isomorphic to
a subquotient of ${\cal F}'_{U'}$.
If one of the following conditions
{\rm (1)} and {\rm (2)} below
is satisfied and if
$f\circ h\colon W\to Y$
is $SS{\cal F}'$-transversal,
then $f\colon X\to Y$
is $SS{\cal F}$-transversal.

{\rm (1)}
The morphism
$h\colon W\to X$
is proper, surjective and generically finite and
the composition
$W\to Y$ is quasi-projective.
The reduced geometric fiber
$U_{\bar k}$ is smooth of dimension $d$
over $\bar k$
and there exists a locally constant sheaf
${\cal G}$ of $\Lambda$-modules
on $U$
such that ${\cal F}_U={\cal G}[d]$.

{\rm (2)}
The morphism
$h$ is quasi-finite
and faithfully flat and $U=X$.
\end{lm}

\proof{
Assume that (1) is satisfied.
Replacing $X$ by the reduced closed subscheme,
we may assume that $X$ is reduced.
Shrinking $U$ if necessary,
we may assume that $h_U\colon U'\to U$
is finite and faithfully flat.
Then, by Corollary \ref{coralt},
${\cal F}$ is isomorphic to a subquotient
of ${}^{\rm p}\! {\cal H}^0Rh_*{\cal F}'$.
Since $W\to Y$ is quasi-projective,
the assertion follows
from Lemma \ref{lmtrfun}.3
and Lemma \ref{lmtr}.3.

Assume that (2) is satisfied.
Since the assertion is \'etale local
on $X$ by Lemma \ref{lmtrb}.1,
we may assume that
$h\colon W\to X$ is finite and faithfully flat
and that $W,X$ and $Y$ are affine.
Then, by Corollary \ref{coralt},
${\cal F}$ is isomorphic to a subquotient
of the perverse sheaf $h_*{\cal F}'$.
Since $W\to Y$ is affine and hence
quasi-projective,
the assertion also follows from
Lemma \ref{lmtrfun}.3 and
Lemma \ref{lmtr}.3.
\qed}

\subsection{Alteration and transversality}\label{ssat}

Let $f\colon X\to Y$
be a morphism of smooth
schemes over a field $k$
and let $D\subset Y$
be a divisor smooth over $k$.
In this article,
we say that $f\colon X\to Y$
is {\em semi-stable}
relatively to $D$
if \'etale locally on $X$
and on $Y$,
there exists a cartesian diagram
$$\begin{CD}
X@>f>>Y@<<<D\\
@VVV
\hspace{-10mm}
\square
\hspace{7mm}
@VVV
\hspace{-10mm}
\square
\hspace{7mm}
@VVV\\
{\mathbf A}^n
@>>>
{\mathbf A}^1@<<<0
\end{CD}$$
where
the lower left horizontal arrow 
${\mathbf A}^n={\rm Spec}\, k[t_1,\ldots,t_n]\to
{\mathbf A}^1
={\rm Spec}\, k[t]$
is defined by $t\mapsto
t_1\cdots t_n$
and the lower right horizontal arrow is
the inclusion of the origin
$0\in {\mathbf A}^1$.
A semi-stable morphism
$f\colon X\to Y$
is flat and 
the base change
$f_V\colon X\times_YV\to V=Y\sm D$
is smooth.
We recall statements on 
the existence of alteration.

\begin{lm}\label{lmdJ}
Let $k$ be a perfect field
and let
$f\colon X\to Y$
be a dominant separated morphism
of integral schemes
of finite type over $k$.

{\rm 1.}
There exists a commutative diagram
\begin{equation}
\begin{CD}
X@<<< W\\
@VfVV@VV{q}V\\
Y@<<<Y'
\end{CD}
\label{eqdJ}
\end{equation}
of integral
schemes of finite type over $k$
satisfying the following condition:
The bottom horizontal arrow
$Y'\to Y$
is dominant
and is the composition $gh$
of an \'etale morphism $g$
and a finite flat radicial morphism $h$.
The schemes $W$ and $Y'$
are smooth over $k$
and the morphism $q\colon W\to Y'$ 
is quasi-projective and smooth.
The induced morphism
$W\to X\times_YY'$
is proper surjective
and generically finite.

{\rm 2.}
Let $\xi\in Y$ be a point
such that the local ring
${\cal O}_{Y,\xi}$
is a discrete valuation ring.
Then, 
there exists a commutative diagram
{\rm (\ref{eqdJ})}
of integral
schemes of finite type over $k$
satisfying the following condition:
The bottom horizontal arrow
$Y'\to Y$ is quasi-finite
and flat and its image is an open neighborhood
of $\xi$.
The schemes $W$ and $Y'$
are smooth over $k$,
the closure $D'\subset Y'$
of the inverse image
of $\xi$ is a divisor smooth over $k$ and
the morphism $g\colon W\to Y'$ 
is quasi-projective
and is semi-stable relatively to $D'$.
The induced morphism
$W\to X\times_YY'$
is proper surjective
and generically finite.
\end{lm}

\proof{
1. Let $\eta$ be the generic point
of $Y$.
Then, it suffices to apply
\cite[Theorem 4.1]{dJ} to
the generic fiber $X\times_Y\eta$.

2.
Let $S={\rm Spec}\
{\cal O}_{Y,\xi}$ be the localization
at $\xi$.
Then, it suffices to apply
\cite[Theorem 8.2]{dJ} to
the base change $X\times_YS\to S$.
\qed}
\medskip

We prove an analogue
of the generic local acyclicity theorem
\cite[Th\'eor\`eme 2.13]{TF}.

\begin{pr}\label{prV1}
Let $f\colon X\to Y$ be
a morphism of schemes 
of finite type over a {\em perfect} field $k$.
Let ${\cal F}$
be a constructible complex 
of $\Lambda$-modules on $X$.
Then, there exists a cartesian diagram
\begin{equation}
\begin{CD}
X@<<< X'\\
@VfVV
\hspace{-10mm}
\square
\hspace{7mm}
@VV{f'}V\\
Y@<<<Y'
\end{CD}
\label{eqdJ2}
\end{equation}
of schemes of finite type over $k$
satisfying the following conditions:
The scheme $Y'$ is smooth over $k$
and the morphism $Y'\to Y$
is dominant and 
is the composition $gh$
of an \'etale morphism $g$
and a finite flat radicial morphism $h$.
The morphism
$f'\colon X'\to Y'$
is properly $SS{\cal F}'$-transversal
for the pull-back
${\cal F}'$ of ${\cal F}$
on $X'$.
\end{pr}

\proof{
We may assume that
${\cal F}$ is a simple
perverse sheaf
by Lemma \ref{lmtr}.3
and Lemma \ref{lmtrb}.
Hence, we may
assume that 
there exist a locally closed immersion
$j\colon Z\to X$ of 
a smooth irreducible 
scheme of dimension $d$
and
a simple locally constant
sheaf ${\cal G}$
of $\Lambda$-modules
such that $j_{!*}{\cal G}[d]={\cal F}$
by \cite[Th\'eor\`eme 4.3.1 (ii)]{BBD}.
By replacing $X$ by the closure
of $j(Z)$,
we may assume that $j\colon
Z\to X$ is an open immersion.
It suffices to consider the
case where 
$Z\to Y$ is dominant
since the assertion is clear
if otherwise.

Let $Z_1\to Z$ be
a finite \'etale covering
such that the pull-back of
${\cal G}$ is constant
and let $X_1$ be the normalization
of $X$ in $Z_1$.
Namely, $X_1$ is the scheme
finite over $X$ corresponding to
the integral closure of
${\cal O}_X$ in 
the quasi-coherent ${\cal O}_X$-algebra
defined as the
direct image of ${\cal O}_{Z_1}$.
The inverse image of $Z\subset X$
in $X_1$ is canonically identified with $Z_1$.
Applying Lemma \ref{lmdJ}.1
to $X_1\to Y$,
we obtain a commutative diagram
\begin{equation}
\begin{CD}
X@<r<<W\\
@VfVV@VVV\\
Y@<<<Y'
\end{CD}
\label{eqWXY}
\end{equation}
of schemes over $k$
satisfying the following conditions:
The scheme $Y'$ is
smooth and the morphism $Y'\to Y$
is dominant and is the composition $gh$
of an \'etale morphism $g$
and a finite flat radicial surjective
morphism $h$.
The morphism
$W\to Y'$ is quasi-projective
and smooth.
The induced morphism
$r'\colon W\to X'=X\times_YY'$ 
is proper surjective and
generically finite.
The pull-back ${\cal G}'_W$ 
of ${\cal G}$ on $W\times_XZ$
is a constant sheaf.

We consider the cartesian diagram
$$\begin{CD}
Z@<<< Z'@<<< W\times_XZ\\
@VjVV 
\hspace{-10mm}
\square
\hspace{7mm}
@V{j'}VV 
\hspace{-15mm}
\square
\hspace{12mm}
@VV{j_W}V\\
X@<<< X'@<{r'}<< W
\end{CD}$$
and let ${\cal G}'$
be the pull-back of ${\cal G}$
on $Z'$.
After shrinking $Z$ if necessary,
we may assume that the reduced
part of $Z'$ is smooth over the perfect field $k$.
Since the finite radicial surjective morphism $h$ is universally
a homeomorphism, we have
${\cal F}'=j'_{!*}{\cal G}'[d]$.

Since 
${\cal G}'_W$
is a constant sheaf on $W\times_XZ$
and $W$ is smooth over $k$,
the intermediate extension
$j_{W!*}{\cal G}'_W[d]$
is constant.
The smooth morphism
$W\to Y'$ is properly
$SSj_{W!*}{\cal G}'_W[d]$-transversal
by Lemma \ref{lmtr}.1.
Since $W\to X'$ is proper
and $W\to Y'$ is quasi-projective,
the morphism $X'\to Y'$ is
$SS{\cal F}'$-transversal by
the case (1) in Lemma \ref{lmtrqf}.
After shrinking $Y'$,
the morphism $X'\to Y'$ is
properly $SS{\cal F}'$-transversal
by Lemma \ref{lmtrV}.2.
\qed}

\begin{cor}\label{corFr}
Let $f\colon X\to Y$ and ${\cal F}$
be as in Proposition {\rm \ref{prV1}}
and assume that $k$
is of characteristic $p>0$.
Then, there exist a
dense open subscheme
$V\subset Y$ smooth over $k$
and an iteration $\tilde V\to V$ of Frobenius
such that the base change
$X\times_Y\tilde V\to \tilde V$
is $SS\tilde {\cal F}$-transversal
for the pull-back $\tilde {\cal F}$
on $\tilde X_V=X\times_Y\tilde V$.
\end{cor}

\proof{
After shrinking $Y'$
in the conclusion of
Proposition \ref{prV1},
we may assume that
$Y'\to Y$
is the composition $jgh$
of an open immersion $j
\colon V\to Y$,
a finite surjective radical morphism
$g$ and an \'etale surjective
morphism $h$.
By Lemma \ref{lmtrb}.1,
we may assume that
$Y'\to Y$ is $jg$.
Thus, the assertion follows.
\qed}

\medskip

We show an analogue of
the stable reduction theorem.

\begin{pr}\label{prV2}
Let 
\begin{equation}
\begin{CD}
X@<\supset<< U\\
@VfVV
\hspace{-10mm}
\square
\hspace{7mm}
@VV{f_V}V\\
Y@<\supset<< V
\end{CD}
\label{eqXUV2}
\end{equation}
be
a cartesian diagram of schemes 
of finite type over a {\em perfect} field $k$.
Assume that $Y$ is normal
and that $V$ is
a dense open subset of $Y$
smooth over $k$.
Let ${\cal F}_U$
be a perverse sheaf 
of $\Lambda$-modules on $U$
such that
$f_V\colon U\to V$
is $SS{\cal F}_U$-transversal.

Then, there exists a cartesian diagram
\begin{equation}
\begin{CD}
X@<<< X'@<{j'}<<U'&=U\times_XX'\\
@VfVV
\hspace{-10mm}
\square
\hspace{7mm}
@VV{f'}V
\hspace{-10mm}
\square
\hspace{7mm}
@VVV\\
Y@<g<<Y'@<\supset<<V'&=V\times_YY'
\end{CD}
\label{eqV2}
\end{equation}
of schemes of finite type over $k$
satisfying the following conditions:
The scheme $Y'$ is smooth over
$k$ and $V'\subset Y'$ is
the complement of
a divisor $D'\subset Y'$
smooth over $k$.
The morphism $g\colon Y'\to Y$
is quasi-finite flat and 
the complement $Y\sm g(Y')$ 
is of codimension $\geqq 2$ in $Y$.
The pull-back 
${\cal F}'_{U'}$ of ${\cal F}_U$ 
is a perverse sheaf on $U'$
and 
for ${\cal F}'=j'_{!*}{\cal F}'_{U'}$ on $X'$,
the morphism
$f'\colon X'\to Y'$
is $SS{\cal F}'$-transversal.
\end{pr}
\medskip

First, we prove a basic case.

\begin{lm}\label{lmst}
Let $X$ and $Y$ be smooth schemes
over a field $k$.
Let $D\subset Y$
be a divisor smooth over $k$
and $V=Y\sm D$ be the complement.
Let $f\colon X\to Y$ be a morphism 
over $k$ semi-stable 
relatively to $D$.
Assume that $\dim X=n$.
For a cartesian diagram {\rm (\ref{eqV2})}
such that $Y'\to Y$
is a quasi-finite flat morphism 
of smooth schemes over $k$, let ${\cal F}'$ be the
perverse sheaf 
${\cal F}'=j'_{!*}\Lambda_{U'}[n]$ on $X'$.

{\rm 1.}
Assume that $\dim Y=1$.
Let $Y'\to Y$
be a flat morphism of smooth curves over $k$
such that for every
$y'\in Y'\sm V'$, 
the action of the
inertia group $I_{y'}$
on $R\Psi_{y'}\Lambda_{U'}$
is trivial.
Then the morphism
$X'\to Y'$ is properly 
$SS{\cal F}'$-transversal.

{\rm 2.}
There exists a quasi-finite faithfully flat morphism
$Y'\to Y$ of smooth schemes over $k$
satisfying the following conditions:
The open subscheme $V'\subset Y'$
is the complement of
a divisor $D'$ smooth over $k$
and the morphism
$X'\to Y'$ is properly 
$SS{\cal F}'$-transversal.
\end{lm}

\proof{
1.
Since the question is \'etale local,
we may assume that
$Y={\mathbf A}^1_k=
{\rm Spec}\ k[t]$,
that
$X=X_n=
{\mathbf A}^n_k={\rm Spec}\ k[t_1,\ldots,t_n]$ 
and that the morphism
$f\colon X\to Y$
is defined by $t\mapsto
t_1\cdots t_n$.
We prove the assertion by induction on $n$.
If $n=1$, then
$f\colon X\to Y$ is \'etale and
${\cal F}'$ is constant.
Hence the assertion follows in this case
by Lemma \ref{lmtr}.1.

Assume $n>1$.
Outside the closed point
$u\in X$ defined by $t_1=\cdots=t_n=0$,
locally there exists a smooth morphism
$X=X_n\to X_{n-1}$ over $Y$.
Hence, the induction hypothesis implies
the assertion on the complement
$X\sm \{u\}$
by Lemma \ref{lmtrfun}.2.
Thus, 
the morphism $f'\colon X'\to Y'$
is properly $SS{\cal F}'$-transversal outside
the inverse image of $u$.
By Proposition \ref{prS}.2
(3)$\Rightarrow$(1),
the morphism $f'\colon X'\to Y'$
is locally acyclic relatively to
${\cal F}'$.
Hence by Lemma \ref{lmtrZ}.2,
the morphism $f'\colon X'\to Y'$
is properly $SS{\cal F}'$-transversal 
on $X'$.

2.
It follows from 1 by Lemma \ref{lmtrb}.1 and Lemma \ref{lmtrV}.1.
\qed}

\proof[Proof of Proposition {\rm \ref{prV2}}]{
The proof is similar to
that of Proposition \ref{prV1}.
By Lemma \ref{lmtrb}
and Lemma \ref{lmtrV},
it suffices to show the
assertion on a neighborhood
of each point $\xi\in Y$
of codimension $1$
not contained in $V$.
Thus, we may assume
that $Y$ is smooth over $k$,
the closure $D$ of
$\xi$ is a divisor smooth over $k$
and that 
$V=Y\sm D$.

We prove the assertion by the induction
on the dimension $d$ of the support
of ${\cal F}_U$.
If $d=0$, it suffices to take $Y'=Y$.
To prove the induction step,
we show the following.

\begin{cl}\label{clind}
{\rm 1.}
Let $Y''\to Y'$ be a morphism of
quasi-finite flat schemes over $Y$
such that $Y'$ and $Y''$ are smooth over
$k$ and that $V'=V\times_YY'$
and $V''=V\times_YY''$
are complements of divisors
$D'$ and $D''$ smooth over $k$.
Assume that ${\cal F}'=j'_{!*}{\cal F}'_{U'}$
satisfies the conclusion of Proposition {\rm \ref{prV2}}.
Then its pull-back ${\cal F}''$
on $X''=X\times_YY''$
also satisfies the same condition.

{\rm 2.}
Let $0\to {\cal G}_U\to
{\cal F}_U\to{\cal H}_U\to0$
be an exact sequence of perverse sheaves
on $U$.
Then, if the assertion of Proposition {\rm \ref{prV2}}
holds for ${\cal G}_U$ and for ${\cal H}_U$,
it also holds for ${\cal F}_U$.
\end{cl}

\proof[Proof of Claim]{
1.
By Lemma \ref{lmtrb},
the $SS{\cal F}'$-transversality
of $X'\to Y'$ 
implies that ${\cal F}''$ is a perverse sheaf on
$X''$
and that
$X''\to Y''$ is
$SS{\cal F}''$-transversal.
By Lemma \ref{lmtrV}.1,
we have
${\cal F}''=j''_{!*}{\cal F}''_{U''}$.

2.
Let $Y'_1$ and
$Y'_2$ be quasi-finite flat schemes
over $Y$ satisfying the conditions in 
Proposition {\rm \ref{prV2}}
for ${\cal G}_U$ and for ${\cal H}_U$
are satisfied respectively.
By Claim 1,
replacing $Y'_1$ and
$Y'_2$ by the normalization $Y'_3$ of
the fiber product
$Y'_1\times_YY'_2$
and shrinking $Y'_3$ if necessary,
we may assume $Y'_1=Y'_2$.
Let $S\to Y'_1$ be the strict localization
at a geometric point $\xi'_1\in Y'_1$
above the generic point $\xi\in D$.
Then, by Proposition \ref{prS}.2,
there exists a finite ramified extension
$S'\to S$ such that 
the pull-back ${\cal F}_{X_{\eta'}}$ of ${\cal F}$
to the generic fiber $X_{\eta'}
\subset X\times_YS'$
satisfies the equivalent conditions loc.~cit.
Hence by Claim 1, further replacing $Y'_1$,
we may assume that $S'=S$.

Let $j_1\colon U'_1=U\times_XX'_1\to X'_1
=X\times_YY'_1$ be the open immersion,
${\cal F}'_{U'_1}$ be the pull-back of 
${\cal F}_U$ on $U'_1$ and let
$A\subset X'_1$
be the union of the supports of constituants
of $j'_{1!*}{\cal F}'_{U'_1}$ that do not meet
$U'_1=X'_1\times_{Y'_1}V'_1$.
Since ${\cal F}_{X_{\eta'}}$ satisfies
the condition (1) in Proposition \ref{prS}.2,
the intersection of $A$ with
the fiber $X'_1\times_{Y'_1}\xi'_1$ is empty.
Hence by replacing $Y'_1$ 
by an open neighborhood $Y'$ of $\xi'_1$,
we may assume that $A$ itself is empty.
Namely, $j'_{!*}{\cal F}'_{U'}$
has no non-trivial subquotient perverse sheaf
supported on $X'\times_{Y'}D'$.

Then, we have an exact sequence
$0\to j'_{!*}{\cal G}'_{U'}
\to j'_{!*}{\cal F}'_{U'}
\to j'_{!*}{\cal H}'_{U'}\to 0$ and
the assertion follows by
Lemma \ref{lmtr}.3.
\qed
}

By Claim and by induction on
the dimension of support,
similarly as
in the proof of Proposition \ref{prV1},
we may
assume that 
there exist a dense open immersion
$j\colon Z\to U$ of 
a smooth irreducible
scheme of dimension $d$
and 
a simple locally constant
sheaf ${\cal G}$
of $\Lambda$-modules
such that ${\cal F}_U=
j_{!*}{\cal G}[d]$.
Further, we may assume
that $Z\to Y$ is dominant.

Taking a finite \'etale
covering trivializing ${\cal G}$
and applying Lemma \ref{lmdJ}.2
as in the proof of
Proposition \ref{prV1},
we obtain a commutative diagram
\begin{equation}
\begin{CD}
X@<r<<W_1\\
@VfVV@VVV\\
Y@<g<<Y_1
\end{CD}
\label{eqWXY2}
\end{equation}
of schemes over $k$
satisfying the following conditions:
The scheme $Y_1$ is
smooth over $k$, the morphism $g_1\colon
Y_1\to Y$
is quasi-finite and flat
and $Y\sm g_1(Y_1)$ is
of codimension $\geqq 2$ in $Y$.
The inverse image $V\times_YY_1$
is the complement of a divisor
$D_1$ smooth over $k$
and the morphism
$W_1\to Y_1$ is quasi-projective
and is semi-stable relatively to $D_1$.
The induced morphism
$r_1\colon W_1\to X_1
=X\times_YY_1$ 
is proper surjective and
generically finite.
The pull-back ${\cal G}'_1$ 
of ${\cal G}$ on $W_1\times_XZ$
is a constant sheaf.

By Lemma \ref{lmst}.2 applied
to the semi-stable morphism 
$W_1\to Y_1$, we obtain
a quasi-finite faithfully flat
morphism $Y'\to Y_1$
of smooth schemes
satisfying the condition loc.\ cit.
We consider the cartesian diagram
$$\begin{CD}
Z@<<< Z'@<<< W'\times_XZ\\
@V{j_Z}VV 
\hspace{-18mm}
\square
\hspace{15mm}
@V{j_{Z'}}VV 
\hspace{-20mm}
\square
\hspace{16mm}
@VV{j_W}V\\
X@<<< X'@<{r'}<< W'\\
@. =X\times_YY'@. =W_1\times_{Y_1}Y'
\end{CD}$$
and let ${\cal G}'$ and
${\cal G}'_{W'}$
denote the pull-backs of
${\cal G}$ on $Z'$
and on $W'\times_XZ$ respectively.
After shrinking $Z$ if necessary,
we may assume that the reduced part of
$Z'$ is smooth over $k$.
Since 
${\cal G}'_{W'}$
is a constant sheaf on $W'\times_XZ$,
the morphism $W'\to Y'$
is $SSj_{W!*}{\cal G}'_{W'}[d]$-transversal
by Lemma \ref{lmst}.2.

The pull-back ${\cal F}'_{U'}$
is a perverse sheaf by Lemma \ref{lmtrb}.
Let ${\cal H}\subset {\cal F}'_{U'}$
be the largest sub perverse sheaf
supported on the complement
$U'\sm Z'$ and let
${\cal H}\subset {\cal H}'\subset {\cal F}'_{U'}$
be the smallest sub perverse sheaf
such that ${\cal F}'_{U'}/{\cal H}'$ 
is supported on $U'\sm Z'$.
Let $u\colon Z'\to U'$ be the open immersion.
Since the restriction of
${\cal F}'_{U'}$ on $Z'$ is identified with
${\cal G}'[d]$,
the subquotient ${\cal H}'/{\cal H}$
is canonically identified
with $u_{!*}{\cal G}'[d]$
by \cite[Corollaire 1.4.25]{BBD}.
Since $r'\colon W'\to X'$
is proper surjective
and generically finite
and since
$r'\colon W'\to Y'$ is quasi-projective,
the morphism $X'\to Y'$ is
$SS{\cal F}'_1$-transversal
for ${\cal F}'_1=j_{Z'!*}{\cal G}'[d]
=j'_{!*}({\cal H}'/{\cal H})$
by the case (1) in Lemma \ref{lmtrqf}.
Since $\dim (U'\sm Z')< d$,
the assertion follows from the
hypothesis of induction and Claim.
\qed}

\begin{cor}\label{corla}
Let the cartesian diagram
{\rm (\ref{eqXUV2})} 
and a perverse sheaf ${\cal F}_U$
on $U=X\times_YV$ be
as in Proposition {\rm \ref{prV2}}.
Assume further that
$Y$ is smooth over $k$ and
that $V$ is the complement
of a divisor $D\subset Y$
smooth over $k$.
Then,
there exists a cartesian diagram
{\rm (\ref{eqV2})} 
satisfying the following conditions:
The scheme $Y'$ is smooth over
$k$ and $V'\subset Y'$ is
the complement of
a divisor $D'\subset Y'$
smooth over $k$.
The morphism $g\colon Y'\to Y$
is quasi-finite flat, 
the morphism
$D'\to D$ is dominant and
the morphism $V'\to V$ is {\em \'etale}.
The pull-back 
${\cal F}'_{U'}$ of ${\cal F}_U$ 
is a perverse sheaf on $U'$
and the morphism
$f'\colon X'\to Y'$
is universally locally acyclic
relatively to ${\cal F}'=j'_{!*}{\cal F}'_{U'}$.
\end{cor}

\proof{
Let $V'\subset Y'$ be
as in the conclusion of
Proposition \ref{prV2}.
Let $\bar Y''$ be
the normalization of $Y$
in the separable closure of
$k(Y)$ in $k(Y')$.
Then, there exists 
a dense open subset
$Y''\subset \bar Y''$ smooth over $k$
of the image of $Y'\to \bar Y''$
such that 
$g''\colon Y''\to Y$ is flat,
that $V''=V\times_YY''$
is the complement of
a divisor $D''$ smooth over $k$,
that $D''\to D$ is dominant,
and that $V''\to V$ is {\em \'etale}.
Since $Y'\times_{\bar Y''}Y''\to
Y''$ is finite surjective
radicial,
the cartesian diagram
(\ref{eqV2}) defined
by $Y''\to Y$ in place of
$Y'\to Y$ satisfies the conditions.
\qed}

\begin{cor}\label{corla1}
Let $X\to Y$ be a morphism
of schemes of finite type
over a field $k$
and assume that 
$Y$ is smooth of {\em dimension} $1$.
Then, for a constructible complex
${\cal F}$ of $\Lambda$-modules
on $X$,
the following conditions are 
equivalent:

{\rm (1)}
$X\to Y$ is locally acyclic
relatively to ${\cal F}$.

{\rm (2)}
$X\to Y$ is universally locally acyclic
relatively to ${\cal F}$.

{\rm (3)}
There exists a finite faithfully
flat morphism 
$Y'\to Y$ of smooth curves
over $k$ such that
the base change
$X'\to Y'$ is $SS{\cal F}'$-transversal
for the base change ${\cal F}'$
of ${\cal F}$ on $X'$.
\end{cor}

The equivalence (1)$\Leftrightarrow$(2) is proved in \cite{Or}.

\proof{
We show (1)$\Rightarrow$(3).
Since the nearby cycles
functor is $t$-exact,
we may assume that ${\cal F}$
is a perverse sheaf.
Then, the assertion 
follows from Propositions \ref{prV1},
\ref{prV2} 
and \ref{prS}.
The implication
(3)$\Rightarrow$(2) is proved in Lemma \ref{lmtr}.2.
The implication
(2)$\Rightarrow$(1) is trivial.
\qed}
\medskip

The following example shows that
taking a covering $Y'\to Y$ in condition (3)
is necessary.

\begin{ex}\label{ex}
{\rm
Let $k$ be a field of characteristic $p\geqq 2$.
Let $X={\mathbf A}^1\times {\mathbf P}^1$
and $j\colon U={\mathbf A}^1\times {\mathbf A}^1
={\rm Spec}\ k[x,y]\to X$
be the open immersion.
Let $W\to X$ be the Artin-Schreier covering
defined by $t^p-t=xy$
ramified along the divisor $X\sm U$
and let 
${\cal G}$ be the locally constant sheaf
of $\Lambda$-modules of
rank $1$ on $U$ trivialized by $W_XU$
and defined by a non-trivial character
${\mathbf F}_p\to \Lambda^\times$.
Let ${\cal F}=j_!{\cal G}$.

Let $W\to X\to Y={\mathbf P}^1$
be the composition with the second projection and
let $Y'\to Y={\mathbf P}^1$
be a finite flat morphism of smooth curves
over $k$ such that the ramification
index at $\infty$ is divisible by $p$.
Then, the normalization $W'\to W\times_YY'$
of the base change is smooth over $Y'$
and $W'\to X'$ is finite and faithfully flat.
Hence by Lemma \ref{lmtrqf}.2,
for the pull-back ${\cal F}'$
of ${\cal F}$
on $X'=X\times_YY'$,
the morphism
$X'\to Y'$ is $SS{\cal F}'$-transversal.
By Corollary \ref{corla1} (3)$\Rightarrow$(2),
$X\to Y$ is local acyclicity relatively to ${\cal F}$ 
(cf.\ \cite[Th\'eor\`eme 2.4.4]{KL}).

On the other hand,
on the complement of the origin
$(0,\infty)\in X$,
the singular support $C=SS{\cal F}$
is the union of the zero-section
$T^*_XX$ and the conormal bundles
$T^*_{X_\infty}X$ of the fiber 
$X_\infty={\rm pr}_2^{-1}(\infty)$.
Hence 
the projection 
${\rm pr}_2\colon X\to Y={\mathbf P}^1$
is not $C$-transversal.
}
\end{ex}

\subsection{Potential transversality}\label{sst}

We prove a refinement of the analogue of
the stable reduction theorem,
using the following consequence
of the stable reduction theorem
for curves.

\begin{lm}\label{lmTem}
Let 
$$
\begin{CD}
U@>{\subset}>>X\\
@V{f_V}VV
\hspace{-10mm}
\square
\hspace{7mm}
@VVfV\\
V@>{\subset}>>
Y@>g>>S
\end{CD}
$$
be a cartesian diagram of morphisms of
smooth schemes of finite type
over a {\em perfect} field $k$
satisfying the following conditions:
The morphism
$f\colon X\to Y$ is flat
and the morphisms
$g\colon Y\to S$ 
and $f_V\colon U\to V$ are 
smooth of relative dimension $1$.
The horizontal arrows are
open immersions and
the open subset $V\subset Y$ is the
complement of 
a divisor $D\subset Y$
smooth over $k$ and
quasi-finite and flat over $S$.

Then, there exists a commutative
diagram
$$
\begin{CD}
X'@>{f'}>>Y'@>{g'}>>S'\\
@VVV@VVV@VVV\\
X@>f>>Y@>g>>S
\end{CD}
$$
of smooth schemes over $k$
satisfying the following conditions:
The morphisms
$S'\to S$, $Y'\to Y\times_SS'$
and $X'\to X\times_YY'$ are 
quasi-finite flat and dominant.
The morphisms $g'\colon Y'\to S'$ 
and $f'\colon X'\to Y'$ are smooth of relative 
dimension $1$
and that $V'=V\times_YY'
\subset Y'$ is the
complement of 
a divisor $D'\subset Y'$
smooth over $k$
and quasi-finite and flat over $S'$.
The morphism
$V'\to V\times_SS'$
is {\em \'etale} and the morphism
$U'=X'\times_{Y'}V' \to U\times_VV'$
is an isomorphism.
The quasi-finite morphisms
$D'\to D$ and
$X'\times_{Y'}D'\to
X\times_YD'$ are dominant.
\end{lm}

\proof{
Let $\bar \eta$
be a geometric point of
$S$ defined by an algebraic
closure of the function 
field of an irreducible
component.
Then, it suffices to apply \cite[Theorem 1.5]{Tem}
to the base change of
$X\to Y\to S$ by $\bar \eta\to S$.
\qed}

\begin{thm}\label{thmV2}
Let 
$$
\begin{CD}
U@>{\subset}>>X\\
@VVV
\hspace{-10mm}
\square
\hspace{7mm}
@VVfV\\
V@>{\subset}>>
Y@>>>S
\end{CD}
$$
be a cartesian
diagram of morphisms of
schemes of finite type
over a {\em perfect} field $k$.
Assume that $Y$ and $S$ are
smooth over $k$,
that 
$Y\to S$ is smooth of relative 
dimension $1$
and that $V\subset Y$ is the
complement of 
a divisor $D$ smooth over $k$ and
quasi-finite and flat over $S$.
Let ${\cal F}_U$ be a 
perverse sheaf 
of $\Lambda$-modules
on $U=X\times_YV$
such that
$U\to V$ is $SS{\cal F}_U$-transversal.

Then, there exists a commutative
diagram
\begin{equation}
\begin{CD}
V'@>\subset>>Y'@>>> S' \\
@VVV\hspace{-10mm}
\square
\hspace{7mm}
@VVV@VVV\\
V@>\subset>>Y@>>> S
\end{CD}
\label{eqYS}
\end{equation}
of smooth schemes over $k$
satisfying the following 
conditions {\rm (1)}
and {\rm (2):}

{\rm (1)}
The morphisms
$S'\to S$ and $Y'\to Y\times_SS'$
are quasi-finite flat and dominant.
The horizontal arrow
$Y'\to S'$ is smooth of relative
dimension $1$.
The left square is cartesian and 
$V'\subset Y'$
is the complement
of a divisor $D'\subset Y'$
smooth over $k$ and 
quasi-finite and flat over $S'$.
The induced morphism
$V'\to V\times_SS'$
is {\em \'etale}
and 
$D'\to D$ is {\em dominant}.

{\rm (2)}
Let 
\begin{equation}
\begin{CD}
U'@>{j'}>>X'@>{f'}>> Y'\\
@VVV
\hspace{-10mm}
\square
\hspace{7mm}
@VVV
\hspace{-10mm}
\square
\hspace{7mm}
@VVV\\
U@>>>X@>f>> Y\\
\end{CD}
\label{eqXX'}
\end{equation}
be a cartesian diagram
and let ${\cal F}'_{U'}$
denote the pull-back of
${\cal F}_U$ on $U'$.
Then the morphism
$f'\colon X'\to Y'$
is $SS{\cal F}'$-transversal
for ${\cal F}'
=j'_{!*}{\cal F}'_{U'}$.
\end{thm}

\proof{
Since the assertion is local on $X$,
we may assume that 
there exists a closed immersion
$i\colon X\to {\mathbf A}^n_Y$
for an integer $n\geqq 0$.
By replacing $X$ and ${\cal F}_U$
by ${\mathbf A}^n_Y$
and ${i|_U}_*{\cal F}_U$
on ${\mathbf A}^n_V$,
we may assume that
$X$ is an open subscheme
of ${\mathbf A}^n_Y$.
We prove the assertion by
induction on $n$.

Assume $n=0$ and hence
$X\to Y$ is an open immersion.
Since the open immersion
$U\to V$ is $SS{\cal F}_U$-transversal,
the singular support
$SS{\cal F}_U$
is a subset of the $0$-section
$T^*_UU$ by Lemma \cite[Lemma 3.6.3]{CC}.
Hence ${\cal F}_U$ is
locally constant by \cite[Lemma 2.1(iii)]{Be}.
Let $U_1\to U$ be a finite \'etale covering
such that the pull-back
of ${\cal F}_U$ is constant.
Let $Y_1$ be the normalization of
$Y$ in $U_1$.
There exists a 
quasi-finite flat and dominant
morphism $S'\to S$
of smooth scheme
such that
the normalization
$Y'$ of $Y_1\times_SS'$
is smooth over $S'$
and that 
$V'\subset Y'$ is the complement
of a divisor $D'$ \'etale over $S'$.
After shrinking $S'$,
we may assume that $Y'\to Y\times_SS'$ is flat.
After shrinking $Y'$ keeping
$D'$ dominant over $D$,
we may assume that
$V'\to V\times_SS'$ \'etale.
Then, the condition (1) is satisfied.
Since ${\cal F}'$ on $X'
\subset Y'$ is constant,
the condition (2) is also satisfied
by Lemma \ref{lmtr}.1.

Assume that $n\geqq 1$
and that the assertion holds
for $n-1$.
For the proof of
the induction step,
we first show the following weaker assertion.

\begin{cl}\label{lmind}
Let 
$X\subset {\mathbf A}^n_Y\to
{\mathbf A}^1_Y$
be a projection
and assume that its restriction 
$U\subset {\mathbf A}^n_V
\to {\mathbf A}^1_V$
is $SS{\cal F}_U$-transversal.
Then, there exist a commutative
diagram {\rm (\ref{eqYS})}
satisfying the condition {\rm (1)}
and 
an open subset
$W'\subset {\mathbf A}^1_{Y'}$
satisfying the following condition:

{\rm (2$'$)}
The intersection
$W'\cap {\mathbf A}^1_{D'}$
is dense in ${\mathbf A}^1_{D'}$.
For the cartesian diagram
{\rm (\ref{eqXX'})}
and for the pull-back ${\cal F}'_{U'}$
of ${\cal F}$ on $U'$ 
and ${\cal F}'=j'_{!*}{\cal F}'_{U'}$
on $X'$,
the morphism
$f'\colon X'\to Y'$
is $SS{\cal F}'$-transversal on 
the inverse image
$X'\times_{{\mathbf A}^1_{Y'}}W'
\subset X'$.
\end{cl}

\proof[Proof of Claim]{
By the induction hypothesis
applied to 
$X\subset {\mathbf A}^n_Y\to {\mathbf A}^1_Y
\to {\mathbf A}^1_S$,
there exists a commutative diagram
$$\begin{CD}
Y_1@>>> S_1\\
@VVV@VVV\\
{\mathbf A}^1_Y 
@>>> {\mathbf A}^1_S
\end{CD}$$
satisfying the conditions
{\rm (1)} and {\rm (2)} in Theorem
\ref{thmV2}.
We consider the cartesian diagram
$$
\begin{CD}
U_1@>{j_1}>>
X_1@>>>Y_1\\
@VVV
\hspace{-10mm}
\square
\hspace{7mm}
@VVV
\hspace{-10mm}
\square
\hspace{7mm}
@VVV\\
U@>>>X
@>>>{\mathbf A}^1_Y
\end{CD}$$
and let ${\cal F}_{U_1}$ 
be the pull-back of ${\cal F}_U$.
Then, 
for ${\cal F}_1=
j_{1!*}{\cal F}_{U_1}$
on $X_1$,
the morphism
$X_1\to Y_1$
is $SS{\cal F}_1$-transversal.
The inverse image
$V_1=V\times_YY_1\subset
Y_1$
is the complement of
a divisor $D_1\subset Y_1$
smooth over $k$
and quasi-finite and flat over $S_1$.
The quasi-finite morphism
$V_1\to V\times_SS_1$ is {\em \'etale} 
and the quasi-finite morphism
$D_1\to {\mathbf A}^1_D$
is dominant.

Since the morphism
$S_1\to {\mathbf A}^1_S$
is quasi-finite and flat,
there exists a
quasi-finite, flat
and dominant morphism
$S'\to S$ 
of smooth schemes over $k$
such that
the normalization
$S'_1$ of $S_1\times_SS'$
is smooth over $S'$
and that
the induced morphism
$S'_1\to S_1$
is also quasi-finite, flat
and dominant.
After shrinking $S'_1$ if necessary,
we may assume that
the morphism
$Y_1\times_{S_1}S'_1
\to 
{\mathbf A}^1_Y\times_{
{\mathbf A}^1_S}S'_1$
of smooth curves over
$S'_1$ is flat.
Hence, by replacing $S,Y,S_1$ and 
$Y_1$ by $S',Y\times_SS',
S'_1$ and $Y_1\times_{S_1}S'_1$,
we may assume that
$S_1\to S$ is smooth
of relative dimension 1.

We consider the commutative
diagram
\begin{equation}
\begin{CD}
V_1@>>> Y_1@>>> S_1\\
@VVV
\hspace{-10mm}
\square
\hspace{7mm}
@VVV@VVV\\
V@>>> Y@>>> S
\end{CD}
\label{eqVW1}
\end{equation}
where the left square is cartesian.
Since $V_1\to
V\times_SS_1$ is \'etale,
the left vertical arrow
$V_1\to V$
is also smooth of relative dimension 1.
The middle vertical arrow
$Y_1\to Y$ is flat.
Hence, by Lemma \ref{lmTem}
applied to (\ref{eqVW1}),
there exists a commutative diagram
$$\begin{CD}
Y'_1@>>>Y'@>>> S'\\
@VVV@VVV@VVV\\
Y_1@>>>Y @>>> S
\end{CD}$$ 
of smooth schemes over $k$
satisfying the following conditions:
The morphisms
$S'\to S$,
$Y'\to Y\times_SS'$
and
$Y'_1\to Y_1\times_YY'$
are quasi-finite flat
and dominant.
The morphisms
$Y'\to S'$
and $Y'_1\to Y'$
are smooth of relative dimension 1.
The inverse image
$V'=V\times_YY'$
is the complement 
$Y'\sm D'$ of
a divisor $D'\subset Y'$ smooth over $k$ and
quasi-finite and flat over $S'$.
The morphism
$V'\to
V\times_SS'$ is {\em \'etale} and
the morphism
$Y'_1\times_{Y'}V'\to
V_1\times_VV'$
is an isomorphism.
The quasi-finite morphisms
$D'\to D$ and $Y'_1\times_{Y'}D'
\to Y_1\times_YD'$
are dominant.
Thus the condition (1)
in Theorem \ref{thmV2}
is satisfied.

The composition
$Y'_1\to Y_1\times_YY'\to
{\mathbf A}^1_{Y'}$
is quasi-finite and flat.
We consider the cartesian diagram
$$\begin{CD}
X'' @>>> Y'_1\\
@VVV
\hspace{-10mm}
\square
\hspace{7mm}
@VVV\\
X'@>>>{\mathbf A}^1_{Y'}
@>>> Y'
\end{CD}$$
and the pull-back
${\cal F}''$ of ${\cal F}_1$
on $X''$.
Then, since
$X_1\to Y_1$ 
is $SS{\cal F}_1$-transversal,
the morphism
$X''\to Y'_1$
is $SS{\cal F}''$-transversal
by Lemma \ref{lmtrb}.1.
Since $Y'_1\to Y'$ is smooth,
the composition
$X''\to Y'$ 
is also $SS{\cal F}''$-transversal
by Lemma \ref{lmtrfun}.1.
Since $Y'_1\to {\mathbf A}^1_{Y'}$
is quasi-finite and flat,
the morphism $f'\colon X'\to Y'$ is 
$SS{\cal F}'$-transversal
on the image of $X''$
by the case (2) in Lemma \ref{lmtrqf}.

Let $W'\subset  {\mathbf A}^1_{Y'}$
be the image of $Y'_1$.
The image of $X''\to X'$ equals
the inverse image
$X'\times_{{\mathbf A}^1_{Y'}}W'$.
Hence,
$f'\colon X'\to Y'$ is 
$SS{\cal F}'$-transversal
on $X'\times_{{\mathbf A}^1_{Y'}}W'
\subset X'$.
Since 
$Y'_1\times_{Y'}D'
\to Y_1\times_YD'$
and $D_1\to {\mathbf A}^1_D$
are dominant,
the intersection
$W'\cap {\mathbf A}^1_{D'}$
is dense in ${\mathbf A}^1_{D'}$.
Thus the condition (2$'$)
in Claim is also satisfied.
\qed}
\smallskip

To complete the proof of
the induction step,
we use the following elementary
lemma.

\begin{lm}\label{lmVC}
Let $X$ be an open subset
of a vector space $V$ 
of dimension $n$ over an infinite field $k$
regarded as a smooth scheme over $k$.
Let $C\subset T^*X$
be a closed conical subset of
dimension $\leqq n$.
Then, there exists an isomorphism
$V\to {\mathbf A}^n$
of vector spaces over $k$
such that the compositions
$X\to V\to {\mathbf A}^n
\to {\mathbf A}^1$
with the projections ${\rm pr}_i,
i=1,\ldots,n$
have at most isolated $C$-characteristic points.
\end{lm}

\proof{
Identify the cotangent bundle
$T^*X$ with the product $X\times V^\vee$
with the dual
and let ${\mathbf P}(C)
\subset 
{\mathbf P}(T^*X)
=
X\times {\mathbf P}(V^\vee)$
be the projectivization.
Then, by the assumption 
$\dim C\leqq n$,
the projection
${\mathbf P}(C)
\to {\mathbf P}(V^\vee)$
is generically finite.
By the assumption that $k$ is infinite,
there exists a basis
$p_1,\ldots,p_n$ of $V^\vee$
such that the fibers
of ${\mathbf P}(C)
\to {\mathbf P}(V^\vee)$
at
$\bar p_1,\ldots,\bar p_n\in
{\mathbf P}(V^\vee)$
are finite.

Since the projectivization $p_i^*\colon 
X\times_{{\mathbf A}^1}T^*{\mathbf A}^1
\to T^*X$ is the section
$X\to
{\mathbf P}(T^*X)
=
X\times {\mathbf P}(V^\vee)$
defined by $\bar p_i$,
the closed subset of $X$
where $p_i$ is not $C$-transversal
equals the image of
$\bar p_i\times_{{\mathbf P}(V^\vee)}
{\mathbf P}(C)\to X$
and is finite
for each $i=1,\ldots,n$.
Hence the product of
$p_1,\ldots,p_n\colon
V\to {\mathbf A}^1$
satisfies the condition.
\qed}
\medskip

Set $C_U=SS{\cal F}_U
\subset T^*U$.
By the assumption that
$U\subset{\mathbf A}^n_V\to V$ is
$SS{\cal F}_U$-transversal,
the morphism 
$T^*U\to T^*U/V$
to the relative cotangent bundle
is finite on $C_U$
by \cite[Lemma 1.2 (ii)]{Be}.
The image $\bar C_U\subset T^*U/V$ 
of $C_U$
and its closure $\bar C\subset T^*X/Y$ 
are closed conical subsets.
Since every irreducible
component of $C_U$
is of dimension $\dim X$,
every irreducible
component of $\bar C_U$
is also of dimension $\dim X$.
Hence, 
for the generic point of
each irreducible component of $Y$, 
the fiber of $\bar C_U$ is
of dimension $\leqq n=\dim X-\dim Y$.
Consequently,
for the generic point of
each irreducible component of
$D\subset Y$, 
the fiber of $\bar C$ is
also of dimension $\leqq n$.

By Lemma \ref{lmVC}
applied to the fibers of
the generic points of 
irreducible components of $D$,
after replacing $S$
by a dense open subset,
there exists a coordinate
of ${\mathbf A}^n_Y\supset X$
such that, 
for each $i=1,\ldots,n$,
there exist
a dense open subset $W_i
\subset {\mathbf A}^1_D$ and
an open neighborhood $X_i
\subset X$
of the inverse image
$W_i\times_{{\mathbf A}^1_Y}X$
by the $i$-th projection ${\rm pr}_i$
satisfying the following
condition:
The inverse image of
$\bar C\subset T^*X/Y$
by the morphism 
$X\times_{{\mathbf A}^1_Y}
T^*{\mathbf A}^1_Y/Y
\to T^*X/Y$
of the relative cotangent bundles
induced by ${\rm pr}_i$
is a subset of the $0$-section on $X_i$.

Then, the restriction
$U\to {\mathbf A}^1_V$ 
of ${\rm pr}_i$ is
$SS{\cal F}_U$-transversal 
on $X_i\cap U$.
By Claim
applied to the restriction 
$X_i\to {\mathbf A}^1_Y$ of
${\rm pr}_i$,
there exist a commutative
diagram {\rm (\ref{eqYS})}
satisfying the condition {\rm (1)}
and for each $i=1,\ldots, n$
a dense open subset
$W'_i\subset {\mathbf A}^1_{Y'}$
satisfying the condition (2$'$).
Hence $X'\to Y'$
is $SS{\cal F}'$-transversal
on the union
$W'=\bigcup_{i=1}^n
{\rm pr}_i^{-1}W'_i
\subset X'
\subset{\mathbf A}^n_{Y'}$
of the inverse images
by the projections.
Since $X'\to Y'$
is $SS{\cal F}'$-transversal
on $U'$,
it is $SS{\cal F}'$-transversal
on $W'\cup U'$.
By shrinking $S'$ if necessary,
we may assume that
$Z'=X'\sm (W'\cup U')
=\prod_{i=1}^n
( {\mathbf A}^1_{D'}
\sm ( {\mathbf A}^1_{D'}\cap W'_i))
\subset {\mathbf A}^n_{D'}$ 
is quasi-finite over $S'$.

By Corollary \ref{corla},
there exists a 
cartesian diagram
$$\begin{CD}
U''@>{j''}>>X''@>{f''}>>Y''@<<<V''\\
@VVV\hspace{-10mm}
\square
\hspace{7mm}
@VVV\hspace{-10mm}
\square
\hspace{7mm}
@VVV\hspace{-10mm}
\square
\hspace{7mm}
@VVV\\
U'@>{j'}>>X'@>{f'}>>Y'@<<<V'
\end{CD}$$
of smooth schemes over $k$
satisfying the following condition:
The morphism
$V''\to V'$ is \'etale
and 
$V''\subset Y''$
is the complement of
a divisor $D''\subset Y''$ smooth over $k$.
The morphism
$D''\to D'$ is dominant.
For the pull-back
${\cal F}''_{U''}$ of ${\cal F}'_{U'}$ 
on $U''$
and ${\cal F}''=j''_{!*}{\cal F}''_{U''}$,
the morphism
$f''\colon X''\to Y''$
is universally locally acyclic relatively
to ${\cal F}''$.

By Lemma \ref{lmtrb}.1
and Lemma \ref{lmtrV}.1,
${\cal F}''$ is the pull-back
of ${\cal F}'$ outside the
inverse image $Z''$ of $Z'$
and $f''\colon X''\to Y''$
is $SS{\cal F}''$-transversal
outside the inverse image
$Z''$.

Let $S''\to S'$ be a quasi-finite
flat dominant morphism
of smooth schemes over $k$
such that
the normalization $Y'''$ of
$Y''\times_{S'}S''$
is smooth over $S''$
of relative dimension
$1$ and
that $V'''=V''\times_{Y''}Y'''$
is the complement of
a divisor $D'''$ smooth over $k$.
Let ${\cal F}'''$
be the pull-back of ${\cal F}''$
on $X'''=X'\times_{Y'}Y'''$.
By Proposition \ref{prV1},
we may assume that
the morphism
$f'''\colon X'''\to S''$
is $SS{\cal F}'''$-transversal.

Since $V''\to V'$ is \'etale, 
the morphism $V''\to S'$ is smooth
and $V'''\to V''\times_{S'}S''$
is an isomorphism.
Hence the morphism
$V'''\to V'\times_{S'}S''$
is \'etale.
The morphism $D'''\to D'$ is dominant.
We consider the commutative diagram
$$\begin{CD}
U'''@>{j'''}>>X'''@>{f'''}>>Y'''@>>>S''\\
@VVV\hspace{-10mm}
\square
\hspace{7mm}
@VVV\hspace{-10mm}
\square
\hspace{7mm}
@VVV@VVV\\
U'@>{j'}>>X'@>{f'}>>Y'@>>>S'
\end{CD}$$
where the left and middle
squares are cartesian.
Then the morphism
$f'''\colon X'''\to Y'''$
is universally locally acyclic
and is $SS{\cal F}'''$-transversal
outside the inverse image
$Z'''$ of $Z''$
quasi-finite over $S''$.

Shrinking $S''$,
we may further assume that
$X'''\to S''$
is properly $SS{\cal F}'''$-transversal
by Lemma \ref{lmtrV}.2.
Then, the morphism
$f'''\colon X'''\to Y'''$ is 
$SS{\cal F}'''$-transversal
by Lemma \ref{lmtrZ}.2.
Further we have an isomorphism
${\cal F}'''
=j'''_{!*}j^{\prime\prime\prime*}
{\cal F}'''$
by Lemma \ref{lmtrV}.1.
Thus, the commutative diagram
$$\begin{CD}
V'''@>\subset>>Y'''@>>>S''\\
@VVV\hspace{-10mm}
\square
\hspace{7mm}
@VVV@VVV\\
V@>\subset>>Y@>>>S
\end{CD}$$
satisfies the required conditions.
\qed}

\begin{cor}\label{corV}
Let $f\colon X\to Y$
be a morphism of scheme
of finite type over a perfect field
$k$.
Assume that $Y$ is smooth of
{\em dimension} $1$.
Let ${\cal F}$ be a constructible
complex of $\Lambda$-modules
on $X$.
Assume that
$f\colon X\to Y$ is 
locally acyclic relatively to ${\cal F}$ and 
that there exists a dense
open subset $V\subset Y$
such that
$f\colon X\to Y$
is $SS{\cal F}$-transversal on $V$.
Then, there exists
a cartesian diagram
$$
\begin{CD}
X@<<< X'\\
@VfVV\hspace{-10mm}
\square
\hspace{7mm}
@VV{f'}V\\
Y@<<< Y'
\end{CD}$$
of morphisms of schemes
of finite type over $k$
satisfying the following condition:
The morphism
$Y'\to Y$ is a
finite generically {\em \'etale} morphism
of smooth curves.
For the pull-back ${\cal F}'$ of 
${\cal F}$ on $X'$, 
the morphism $f'\colon X'\to Y'$
is $SS{\cal F}'$-transversal.
\end{cor}

\proof{
Since the shifted vanishing cycles functor
$R\Phi[-1]$ is $t$-exact,
we may assume that
${\cal F}$ is a perverse sheaf.
Then the assertion follows
from the case
where $S={\rm Spec}\ k$
in Theorem \ref{thmV2},
Proposition \ref{prS}.1
and Lemma \ref{lmwa}.
\qed}

\section{Characteristic cycles and the direct image}\label{spp}

\subsection{Direct image
of a cycle}\label{ssdc}

To state the compatibility with
push-forward, we fix some terminology
and notations.
Recall that a morphism
$f\colon X\to Y$ of
noetherian schemes is said to be proper
on a closed subset $Z\subset X$
if its restriction $Z\to Y$
is proper with respect to a
closed subscheme structure of $Z\subset X$.

Let $f\colon X\to Y$ be 
a morphism of smooth schemes over 
a field $k$ and we
consider the diagram
\begin{equation}
\begin{CD}
T^*X@<<< X\times_YT^*Y
@>>> T^*Y
\end{CD}
\label{eqJD}
\end{equation}
as an algebraic correspondence
from $T^*X$ to $T^*Y$.
Assume that
every irreducible component
of $X$ (resp.\ of $Y$)
is of dimension $n$ (resp.\ $m$).
Let $B\subset X$
be a closed subset
on which $f\colon X\to Y$
is proper
and let $C\subset T^*X$
be a closed subset of 
$B\times_XT^*X$.
Then, the closed subset
$f_\circ C\subset
T^*Y$
is defined as the image
by the right arrow
in (\ref{eqJD})
of the inverse image of
$C$ by the left arrow.
It is a closed subset by
the assumption that
$f$ is proper on $B$.
The composition of
the Gysin map 
\cite[6.6]{Ful}
for the first arrow and
the push-forward map
for the second arrow
defines a morphism
\begin{equation}
\begin{CD}
f_!\colon {\rm CH}_n(C)
@>>>
{\rm CH}_m(f_\circ C)
\end{CD}
\label{eqCH}
\end{equation}
since $\dim T^*X
-\dim X\times_YT^*Y
=n-m$.
If every irreducible component
of $C$ 
(resp.\ $f_\circ C$) is of dimension 
$\leqq n$
(resp.\ $\leqq m$),
the morphism (\ref{eqCH})
defines
a morphism
\begin{equation}
\begin{CD}
f_!\colon Z_n(C)
@>>>
Z_m(f_\circ C)
\end{CD}
\end{equation}
of free abelian groups of cycles.

\begin{lm}\label{lmCHtr}
Let $$
\xymatrix{
X\ar[r]^g
\ar[rd]_f &
X'\ar[d]^{f'}\\
&Y}$$ be 
a commutative diagram of
morphisms of smooth schemes over $k$.
Assume that
every irreducible component
of $X$ (resp.\ of $X'$ and $Y$)
is of dimension $n$ (resp.\ 
$n'$ and $m$).
Let $B\subset X$
be a closed subset
on which $f\colon X\to Y$
is proper
and let $C\subset T^*X$
be a closed subset of 
$B\times_XT^*X$.
Then, the diagram
$$
\xymatrix{
{\rm CH}_n(C)\ar[r]^{g_!}
\ar[rd]_{f_!} &
{\rm CH}_{n'}(g_\circ C)
\ar[d]^{f'_!}\\
&{\rm CH}_m(f_\circ C)}$$
is commutative.
\end{lm}

\proof{
We consider the diagram
\begin{equation*}
\begin{CD}
T^*X@<<< X\times_{X'}T^*X'
@>>> T^*X'@.\\
@. @AAA
\hspace{-20mm}
\square
\hspace{17mm}
@AAA\\
@.X\times_YT^*Y
@>>>X'\times_YT^*Y
@>>> T^*Y
\end{CD}
\end{equation*}
with cartesian square.
After decomposing the
right vertical arrow into the composition
of a smooth morphism and a regular immersion,
it suffices to apply \cite[Theorem 6.2 (a)]{Ful}.
\qed}

\begin{lm}\label{lmptrA}
Let $f\colon X\to Y$
be a smooth morphism of smooth irreducible schemes
over a perfect field $k$.
Assume that $X$ 
(resp.\ of $Y$) is
of dimension $n$ 
(resp.\ $m$).
Let $C=\bigcup_aC_a\subset T^*X$
be a closed conical subset such
that every irreducible component
$C_a$ is of dimension $n$
and that $f\colon X\to Y$ is
properly $C$-transversal and
is proper on the base $B=C\cap T^*_XX\subset X$.

Let $A=\sum_am_aC_a$ be a linear combination.
Let $y\in Y$ be a closed point,
let $A_y=i_y^!A$ be the pull-back 
{\rm \cite[Definition 7.1]{CC}} by
the closed immersion $i_y\colon X_y\to X$
of the fiber
and let
$(A_y,T^*_{X_y}X_y)_{T^*X_y}$
denote the intersection number.
Then, we have
\begin{equation}
f_! A=
(-1)^m
(A_y,T^*_{X_y}X_y)_{T^*X_y}
\cdot
[T^*_YY]
\label{eqAy0}
\end{equation}
in $Z_m(T^*_YY)$.
\end{lm}

\proof{
Since the closed immersion
$i_y\colon X_y\to X$
is properly $C$-transversal by Lemma \ref{lmtrbc}.1,
the pull-back
$A_y=i_y^!A$
is defined.
Further by the assumption that
$f\colon X\to Y$
is $C$-transversal,
we have an inclusion
$f_\circ C\subset T^*_YY$
and $f_!A$
is defined as an element of
${\rm CH}_m(f_\circ C)
=Z_m(f_\circ C)\subset Z_m(T^*_YY)$.
Hence it suffices to 
show that the
coefficient of $T^*_YY$
in $f_!A$ equals the intersection number
$(-1)^m(A_y,T^*_{X_y}X_y)_{T^*X_y}$.

We consider the cartesian diagram
$$\begin{CD}
T^*X@<<< X\times_YT^*Y@>>> T^*Y\\
@AAA
\hspace{-20mm}
\square
\hspace{15mm}
@AAA
\hspace{-20mm}
\square
\hspace{15mm}
@AAA\\
X_y\times_XT^*X@<<< X_y\times_YT^*Y@>>> y\times_YT^*Y\\
@VVV
\hspace{-20mm}
\square
\hspace{15mm}
@VVV
\hspace{-20mm}
\square
\hspace{15mm}
@VVV\\
T^*X_y@<0<< X_y@>>> y.
\end{CD}$$
We regard the four sides of the exterior square of the diagram
as algebraic correspondences.
Since $f$ is assumed properly $C$-transversal,
$f_!A$ on the upper right corner
is some multiple of the $0$-section $T^*_YY$.
Hence, the coefficient of $f_!A$
is the image of $A$ by the composition
via the upper right corner.
It equals the composition
via the lower left corner
by \cite[Theorem 6.2 (a)]{Ful}
applied to the upper right and the lower left squares.
The image of $A$ on the lower left corner
is $(-1)^m$-times $i_y^!A$
since the definition of $i_y^!A$ in
\cite[Definition 7.1]{CC} involves
the sign $(-1)^{\dim X-\dim X_y}=(-1)^m$.
Since the bottom line defines
the intersection number 
$(-,T^*_{X_y}X_y)_{T^*X_y}$,
the assertion follows.
\qed}
\medskip

We study the case where $Y$
is a smooth curve
and $\dim f_\circ C=1$.
Let $f\colon X\to Y$
be a morphism of smooth schemes
over $k$.
Assume that
every irreducible
component of $X$ 
(resp.\ of $Y$) is
of dimension $n$ 
(resp.\ $1$).
Let $C\subset T^*X$
be a closed conical subset
such that every irreducible
component $C_a$ of 
$C=\bigcup_aC_a$
is of dimension $n$ and that $f\colon X\to Y$
is proper on the base
$B=C\cap T^*_XX\subset X$.
Let $V\subset Y$ be a dense open subscheme
such that the base change
$f_V\colon X_V\to V$
is properly $C_V$-transversal
for the restriction 
$C_V$ of $C$ on $X_V$.

Let $y\in Y\sm V$ be a closed point
on the boundary and
let $t$ be a uniformizer at $y$
and let $df$ denote the section
of $T^*X$ defined on a neighborhood
of the fiber $X_y$ by the pull-back $f^*dt$.
Then, on a neighborhood of
$X_y$, 
the intersection $C\cap df
\subset T^*X$ is supported on the 
inverse image of the intersection
$B\cap X_y$.
Hence for a linear combination
$A=\sum_am_aC_a$,
the intersection product
\begin{equation}
(A,df)_{T^*X,X_y}
\label{eqAdf}
\end{equation}
supported on the fiber $X_y$
is defined as an element of
${\rm CH}_0(B\cap X_y)$.
Since $C$ is conical,
the intersection product
$(A,df)_{T^*X,X_y}$
does not depend on the choice of $t$.
Thus the intersection number
also denoted $(A,df)_{T^*X,X_y}$
is defined as its image by
the degree mapping
${\rm CH}_0(B\cap X_y)
\to
{\rm CH}_0(y)={\mathbf Z}$.

\begin{lm}\label{lmd1A}
Let $f\colon X\to Y$
be a morphism of smooth irreducible schemes
over a perfect field $k$.
Assume that $X$ 
(resp.\ of $Y$) is
of dimension $n$ 
(resp.\ $1$).
Let $C=\bigcup_aC_a\subset T^*X$
be a closed conical subset
as in Lemma {\rm \ref{lmptrA}}.

{\rm 1.}
The following conditions
are equivalent:

{\rm (1)}
$\dim f_\circ C\leqq 1$.

{\rm (2)}
There exists a dense open subscheme
$V\subset Y$
such that the base change
$f_V\colon X_V\to V$
is $C_V$-transversal
for the restriction 
$C_V$ of $C$ on $X_V$.

{\rm (3)}
There exists a dense open subscheme
$V\subset Y$
such that the base change
$f_V\colon X_V\to V$
is properly $C_V$-transversal
for the restriction 
$C_V$ of $C$ on $X_V$.

{\rm 2.}
Let $V\subset Y$ be a dense open subscheme
satisfying the condition {\rm (3)}
above.
Let $A=\sum_am_aC_a$
be a linear combination,
let $v\in V$ be a closed point
and define the intersection number
$(A_v,T^*_{X_v}X_v)_{T^*X_v}$
as in Lemma {\rm \ref{lmptrA}}.
Then, we have
\begin{equation}
f_! A=
-
(A_v,T^*_{X_v}X_v)_{T^*X_v}\cdot
[T^*_YY]+
\sum_{y\in Y\sm V}
(A,df)_{T^*X,X_y}\cdot [T^*_yY]
\label{eqAy}
\end{equation}
in $Z_1(f_\circ C)$.
\end{lm}

\proof{
1.
Since $f_\circ C$ is a closed conical subset
of the line bundle $T^*Y$,
the condition (1) is equivalent
to the existence of a dense open subset
$V\subset Y$ such that
$f_\circ C\subset T^*_YY
\cup\bigcup_{y\in Y\sm V}T^*_yY$.
This is equivalent to the condition (2).
The equivalence
(2)$\Leftrightarrow$(3) follows
from Lemma \ref{lmtrV}.2.

2.
It suffices to compare the coefficients
of the $0$-section
$T^*_YY$
and of the fibers
$T^*_yY$ respectively.
For those of $T^*_YY$,
it is proved in Lemma \ref{lmptrA}.
For those of $T^*_yY$,
it follows from
the projection formula
\cite[Theorem 6.2 (a)]{Ful}
applied to
the cartesian square in the 
commutative diagram
$$
\xymatrix{
T^*X&& X\times_YT^*Y
\ar[rr]\ar[ll]&&
T^*Y\\
&&X\ar[u]_{\hspace{15mm}\square}\ar[ull]^{df}\ar[rr]
&&Y\ar[u]_{dt}.
}$$
\qed}
\medskip

\begin{lm}\label{lmcE}
Let $X$ be a scheme of
finite type of dimension $d$
over a field $k$
and let $E$ be a vector bundle
on $X$ associated
to a locally free ${\cal O}_X$-module
${\cal E}$ of rank $n$.
Let $s\colon X\to E$
be a section,
$0\colon X\to E$ be the zero section
and
$Z=Z(s)
=0(X)\cap s(X)\subset X$ be the zero locus
of $s$.
Let ${\cal K}=[{\cal O}_X
\overset s\to {\cal E}]$
be the complex of ${\cal O}_X$-modules
where ${\cal E}$ is put on
degree $0$ and let
${c_n}^X_Z({\cal K})$ be the localized Chern
class defined in {\rm \cite[Section 1]{Bl}}.
Then, we have
$$(0(X),s(X))_E
=
{c_n}^X_Z({\cal K})\cap [X]$$
in ${\rm CH}_{d-n}(Z)$.
\end{lm}

\proof{
We may assume that $X$ is integral
and $Z\subsetneqq X$.
By taking the blow-up at $Z$
and by \cite[Proposition 2.3.1.6]{KS},
we may assume that $Z$ is a Cartier divisor $D
\subset X$.
Then, we have an exact sequence of
$0\to {\cal L}\to {\cal E}\to {\cal F}\to 0$
of locally free ${\cal O}_X$-modules
where ${\cal L}$ and ${\cal F}$ are
of rank 1 and $n-1$ respectively
and $s\in \Gamma(X,{\cal E})$
is defined by
$s\in \Gamma(X,{\cal L})$.
Then, the right hand side
equals $c_{n-1}({\cal F})\cap [D]$
by \cite[Proposition (1.1) (iii)]{Bl}.
The left hand side also equals
$c_{n-1}({\cal F})\cap [D]$
by the excess intersection formula
\cite[Theorem 6.3]{Ful}.
\qed}
\medskip

We define the specialization of
a cycle.
Let $f\colon X\to Y$
be a smooth morphism
of smooth schemes over 
a perfect field $k$ and 
assume that $X$ (resp.\ $Y$)
is equidimensional of dimension
$n+1$ (resp.\ 1).
Let $y\in Y$ be a closed point,
$V=Y\sm \{y\}$
be the complement
and $U=X\times_YV$ 
be the inverse image.
Let $C\subset T^*U$
be a closed conical subset
equidimensional of dimension
$n+1$
and assume that
$U\to V$ is properly $C$-transversal.
We define its specialization
$${\rm sp}_yC
\subset T^*X_y$$
as follows.
By the assumption that
$U\to V$ is properly $C$-transversal
and \cite[Lemma 3.1]{CC},
the morphism
$T^*U\to T^*U/V$
to the relative cotangent bundle
is finite on $C$.
Hence its image
$C'\subset T^*U/V$
is a closed conical subset.
Let $\bar C'\subset T^*X/Y$
be the closure and
define
${\rm sp}_yC
\subset T^*X_y$
to be the fiber 
$\bar C'\times_Yy
\subset T^*X/Y\times_Yy
=T^*X_y$.
The specialization
${\rm sp}_yC\subset T^*X_y$
is a closed conical subset
equidimensional of
dimension $n$.

For a linear combination
$A=\sum_am_aC_a$
of irreducible components
of $C=\bigcup_aC_a$,
we define its specialization
$${\rm sp}_yA
\in Z_n({\rm sp}_yC)$$
as follows.
First, we define
$A'\in Z_{n+1}(C')$
as the push-forward of
$A$ by the morphism
$T^*U\to T^*U/V$ finite on $C$.
Let $\bar A'\in Z_{n+1}(\bar C')$
be the unique element
extending $A'\in Z_{n+1}(C')$.
Then, we define
${\rm sp}_yA
\in Z_n({\rm sp}_yC)$
to be the {\em minus} of
the pull-back of $A'$
by the Gysin map
for the immersion
$i_y\colon X_y\to X$.
If $X\to Y$ is proper,
for a closed point $v\in V$
and the closed immersion
$i_v\colon X_v\to X$,
we have
\begin{equation}
({\rm sp}_yA,
T^*_{X_y}X_y)_{T^*X_y}
=
(i_v^!A,
T^*_{X_v}X_v)_{T^*X_v}
\label{eqspiv}
\end{equation}
since the definition of $i_v^!A$ in 
\cite[Definition 7.1]{CC} involves
the sign $(-1)^{\dim X-\dim X_v}=-1$. 

\begin{lm}\label{lmsp}
Let $f\colon X\to Y$
be a smooth morphism
of smooth schemes over 
a perfect field $k$ and 
assume that $X$ (resp.\ $Y$)
is equidimensional of dimension
$n+1$ (resp.\ $1$).
Let $y\in Y$ be a closed point,
$i_y\colon X_y\to X$
be the closed immersion 
of the fiber,
$V=Y\sm \{y\}$
be the complement
and $U=X\times_YV$ 
be the inverse image.
Let $C\subset T^*X$
be a closed conical subset
equidimensional of dimension
$n+1$ such that
$f\colon X\to Y$ is properly $C$-transversal.

{\rm 1.}
For the restriction $C_U$
of $C$ on $U$, we have
\begin{equation}
{\rm sp}_yC_U=
i_y^\circ C.
\label{eqspC}
\end{equation}

{\rm 2.}
For a linear combination
$A=\sum_am_aC_a$
of irreducible components
of $C=\bigcup_aC_a$
and its restriction $A_U$ on $U$,
we have
\begin{equation}
{\rm sp}_yA_U=
i_y^!A.
\label{eqspA}
\end{equation}
\end{lm}

\proof{
1.
By the assumption that
$f\colon X\to Y$ is properly
$C$-transversal
and \cite[Lemma 3.1]{CC},
the morphism
$T^*X\to T^*X/Y$
to the relative cotangent bundle
is finite on $C$
and hence its image
$C'\subset T^*X/Y$
is a closed conical subset.
Further $C'$ with reduced scheme
structure is flat over $Y$.
Hence it equals the closure
of the restriction $C'_U$
and we obtain (\ref{eqspC}).

2.
We consider the cartesian diagram
$$\begin{CD}
T^*X@>>> T^*X/Y\\
@AAA\hspace{-17mm}
\square\hspace{12mm} @AAA\\
X_y\times_XT^*X@>>> T^*X_y.
\end{CD}$$
The right hand side is
the minus of the image of
$A$ by the push-forward
and the pull-back via
upper right.
The left hand side is
the minus of the image of
$A$ by the pull-back and
the push forward
via lower left.
Hence the assertion follows from
the projection formula
\cite[Theorem 6.2 (a)]{Ful}.
\qed}

\subsection{Characteristic cycle of
the direct image}\label{ssccd}

Let $k$ be a field and 
let $\Lambda$ be a finite field
of characteristic $\ell$ invertible in $k$.
Let $X$ be a smooth scheme over $k$
such that every irreducible
component is of dimension $n$.
Let ${\cal F}$ be a constructible complex of $\Lambda$-modules 
on $X$
and $C=SS{\cal F}$ be the singular support.
Then, every irreducible
component $C_a$ of $C=
\bigcup_aC_a$
has the same dimension 
as $X$
\cite[Theorem 1.3 (ii)]{Be} and
the base $B=C\cap T^*_XX
\subset T^*_XX=X$
defined as the intersection
with the $0$-section
equals the support of ${\cal F}$
\cite[Lemma 2.1 (i)]{Be}.
Let $f\colon X\to Y$
be a morphism of smooth schemes
over $k$, proper on the
support of ${\cal F}$.
Then, we have an inclusion
\begin{equation}
SSRf_*{\cal F}\subset f_\circ SS{\cal F}
\end{equation}
by 
\cite[Lemma 2.2 (ii)]{Be}.

We restate a conjecture
from \cite[Conjecture 1]{prop}.

\begin{cn}\label{cnf*}
Let $f\colon X\to Y$
be a morphism 
of smooth schemes 
over a perfect field $k$.
Assume that
every irreducible
component of $X$ 
(resp.\ of $Y$) is
of dimension $n$ 
(resp.\ $m$).
Let ${\cal F}$ be a constructible complex of $\Lambda$-modules 
on $X$
and $C=SS{\cal F}$ be the singular support.
Assume that $f$ is proper on
the support of ${\cal F}$.
Then, we have
\begin{equation}
CCRf_*{\cal F}=
f_!CC{\cal F}
\label{eqcnf}
\end{equation}
in ${\rm CH}_m(f_\circ SS{\cal F})$.
\end{cn}

If $\dim f_\circ SS{\cal F}\leqq m$,
the equality (\ref{eqcnf})
is an equality as cycles
in ${\rm CH}_m(f_\circ SS{\cal F})
=Z_m(f_\circ SS{\cal F})$
without rational equivalence.

A weaker version of Conjecture \ref{cnf*}
is proved in the case $k$ is finite
and $X$ and $Y$ are projective
in \cite{UYZ}
using $\varepsilon$-factors.

If $Y={\rm Spec}\ k$,
the equality (\ref{eqcnf}) means
the index formula
\begin{equation}
\chi(X_{\bar k},{\cal F})
= (CC{\cal F},T^*_XX)_{T^*X}
\label{eqind}
\end{equation}
where the right hand side denotes
the intersection number.
Further if $X$ is projective,
the equality (\ref{eqind})
is proved in
\cite[Theorem 7.13]{CC}.

\begin{lm}\label{lmfun}
Let $f\colon X\to Y$
be a morphism of smooth schemes
over $k$
and let ${\cal F}$
be a constructible
complex of $\Lambda$-modules.
Assume that $f\colon X\to Y$
is proper on the
support of ${\cal F}$.
Assume that
every irreducible
component of $X$ 
(resp.\ of $Y$) is
of dimension $n$ 
(resp.\ $m$).

{\rm 1.}
Let
$$
\xymatrix{
X\ar[r]^g
\ar[rd]_f &
X'\ar[d]^{f'}\\
&Y}$$ be 
a commutative diagram of
morphisms of smooth schemes over $k$.
Then, we have
\begin{equation}
f_!CC{\cal F}=
f'_!(g_!CC{\cal F})
\label{eqfg}
\end{equation}
in ${\rm CH}_m(f_\circ SS{\cal F})$.

{\rm 2.}
Assume that one of the following conditions
{\rm (1)} and {\rm (2)} is 
satisfied:

{\rm (1)} $f\colon X\to Y$
is an immersion.

{\rm (2)} 
$f\colon X\to Y$ is
quasi-projective and
$SS{\cal F}$-transversal.

\noindent
Then, 
we have $\dim (f_\circ SS{\cal F})
\leqq m=\dim Y$
and 
\begin{equation*}
CCRf_*{\cal F}=
f_!CC{\cal F}
\leqno{(\ref{eqcnf})}
\end{equation*}
in $Z_m(f_\circ SS{\cal F})$.
\end{lm}

\proof{
1.
It follows from
Lemma \ref{lmCHtr}.

2.
The case (1) is proved in
\cite[Lemma 5.13.2]{CC}.
We show the case (2).
Since $f_\circ SS{\cal F}$
is a subset of the $0$-section $T^*_YY$,
we have $\dim (f_\circ SS{\cal F})
\leqq m=\dim Y$.
We may assume that $Y$ is connected and affine
and hence $X$ is quasi-projective.
Let $X\to P$ be an immersion
to a projective space
and factorize $X\to Y$
as the composition of an immersion
$X\to Y\times P$
and the projection $Y\times P\to Y$.
Then, by 1 and the case (1),
we may assume
that $f\colon X\to Y$ is projective and smooth.

By the assumption
that $f\colon X\to Y$ is
$SS{\cal F}$-transversal,
it is locally acyclic relatively
to ${\cal F}$ by Lemma \ref{lmtr}.1.
Since $f\colon X\to Y$ is proper,
the direct image
$Rf_*{\cal F}$ is locally constant 
by \cite[5.2.4]{app}.
By Lemma \ref{lmtrV}.2,
there exists a dense
open subscheme $V\subset Y$
such that
$f\colon X\to Y$
is properly $SS{\cal F}$-transversal
on $V$.
By \cite[Lemma 5.11.1]{CC}
and Lemma \ref{lmptrA}, 
it suffices to show
the equality 
\begin{equation}
{\rm rank}\ Rf_*{\cal F}
= (i_y^!CC{\cal F},T^*_{X_y}X_y)
_{T^*X_y}
\label{eqptr}
\end{equation}
for a closed point $y\in V$.
Since 
${\rm rank}\ Rf_*{\cal F}
=\chi(X_{\bar y},i_y^*{\cal F})$,
the equality
(\ref{eqptr}) follows from the compatibility
$CCi_y^*{\cal F}=
i_y^!CC{\cal F}$
with the pull-back \cite[Theorem 7.6]{CC}
and 
the index formula \cite[Theorem 7.13]{CC}.
\qed}
\medskip

We consider the case where $Y$
is a smooth curve
and $\dim f_\circ SS{\cal F}\leqq 1$.
We recall the definition
of the Artin conductor
and the description of
the characteristic cycle
of a sheaf on a curve.
Let $Y$ be a smooth irreducible curve over
a perfect field $k$
and let ${\cal G}$ be
a constructible complex of
$\Lambda$-modules on $Y$.
Let $V\subset Y$ be a dense open subscheme
such that the restriction 
${\cal G}_V$ is locally constant
i.\ e.\
the cohomology sheaf
${\cal H}^q{\cal G}_V$ is locally constant
for every integer $q$.
For a closed point $y\in Y$,
the Artin conductor
$a_y{\cal G}$ is defined by
\begin{equation}
a_y{\cal G}
=
{\rm rank}\ 
{\cal G}_V
-{\rm rank}\
{\cal G}_{\bar y}
+
{\rm Sw}_y{\cal G}.
\label{eqay}
\end{equation}
Here
$\bar y$ denotes a geometric point above $y$
and ${\rm Sw}_y$
denotes the alternating sum of the Swan conductor.
The characteristic cycle is given by
\begin{equation}
CC{\cal G}
=-\Bigr({\rm rank}\ {\cal G}_V
\cdot [T^*_YY]
+\sum_{y\in Y\sm V}a_y{\cal G}\cdot
[T^*_yY]\Bigr)
\label{eqccY}
\end{equation}
by \cite[Lemma 5.11.3]{CC}.
Here 
$T^*_yY$ is the fiber of $T^*Y$ at $y$.

Let $f\colon X\to Y$
be a morphism
of smooth schemes
over a perfect field $k$
and $y\in Y$ be a closed point.
Assume that $\dim Y=1$.
Let ${\cal F}$
be a constructible complex
of $\Lambda$-modules on $X$.
Assume that
$f\colon X\to Y$ is proper on
the support of ${\cal F}$.
Under the assumption
$\dim f_\circ SS{\cal F}\leqq 1$,
the equality
$CCRf_*{\cal F}=
f_!CC{\cal F}$
(\ref{eqcnf})
in $Z_1(f_\circ SS{\cal F})$
is equivalent to
the equality
\begin{equation}
-a_yRf_*{\cal F}=
(CC{\cal F},df)_{T^*X,X_y}.
\label{eqcf}
\end{equation}
for every closed point
$y\subset Y\sm V$
by (\ref{eqccY}),
Lemma \ref{lmfun}.2 (2)
and Lemma \ref{lmd1A}.2,
where the right hand side is defined
as in (\ref{eqAdf}).

\begin{thm}\label{thmay}
Let $f\colon X\to Y$
be a quasi-projective morphism
of smooth schemes
over a perfect field $k$
and $y\in Y$ be a closed point.
Assume that $\dim Y=1$.
Let ${\cal F}$
be a constructible complex
of $\Lambda$-modules on $X$.
Assume that
$f\colon X\to Y$ is proper on
the support of ${\cal F}$
and is properly $SS{\cal F}$-transversal
on a dense open subscheme
$V\subset Y$.
Then, we have
\begin{equation*}
-a_yRf_*{\cal F}=
(CC{\cal F},df)_{T^*X,X_y}.
\leqno{(\ref{eqcf})}
\end{equation*}
\end{thm}

\proof{
We may assume that $k$
is algebraically closed.
By the same argument as in the proof
of Lemma \ref{lmfun}.2,
we may assume that
$f\colon X=Y\times P\to Y$
is the projection
for a projective space $P$.
By Lemma \ref{lmd1A}.1
and by replacing $Y$
by a projective smooth
curve over $k$
containing $Y$ as a dense
open subscheme,
we may assume that $Y$
is projective
and smooth.

By Lemma \ref{lmfun}
applied to $X\to Y\to {\rm Spec}\ k$,
we obtain
$$(f_!CC{\cal F},T^*_YY)_{T^*Y}
=
(CC{\cal F},T^*_XX)_{T^*X}.$$
By the index formula \cite[Theorem 7.13]{CC},
we have
$$
(CCRf_*{\cal F},T^*_YY)_{T^*Y}
=\chi(Y_{\bar k},Rf_*{\cal F})
=\chi(X_{\bar k},{\cal F})
=(CC{\cal F},T^*_XX)_{T^*X}.
$$
Thus, we have
$$
(CCRf_*{\cal F}-
f_!CC{\cal F},T^*_YY)_{T^*Y}
=0.$$
Since the coefficients of
$T^*_YY$ in
$CCRf_*{\cal F}$
and 
$f_!CC{\cal F}$ are equal
by (\ref{eqAy0}), (\ref{eqccY})
and the index formula \cite[Theorem 7.13]{CC},
we obtain
\begin{equation}
\sum_{y\in Y\sm V}
-a_yRf_*{\cal F}=
\sum_{y\in Y\sm V}
(CC{\cal F},df)_{T^*X,X_y}.
\label{eqsum}
\end{equation}

By d\'evissage using 
Lemma \ref{lmtr}.3
and \cite[Lemma 5.13.1]{CC}, we may assume
that ${\cal F}$ is a perverse sheaf.
Set $\widetilde V=V\cup \{y\}$
and $Z=Y\sm \widetilde V$.
By Corollary \ref{corS}, Corollary \ref{corV}
and Lemma \ref{lmwa},
there exists a faithfully flat
finite morphism
$Y'\to Y$ of projective smooth curves
{\em \'etale at} $y$
satisfying the following condition:
Let 
$$\begin{CD}
X@<<< X'@<{\tilde j'}<<
X'_{\widetilde V'}\\
@VfVV\hspace{-10mm}
\square
\hspace{6mm}
@V{f'}VV
\hspace{-13mm}
\square
\hspace{8mm}
@VVV\\
Y@<<< Y'@<<<\widetilde V'&=
Y'\times_Y\widetilde V
\end{CD}$$
be a cartesian diagram
and set
${\cal F}'=\tilde j'_{!*}
{\cal F}'_{\widetilde V'}$
for the pull-back
${\cal F}'_{\widetilde V'}$
of ${\cal F}$ on $X'_{\widetilde V'}$.
Then on $Y'_0=Y'\sm Y'\times_Yy$, the morphism
$f'\colon X'\to Y'$
is $SS{\cal F}'$-transversal and hence is
universally locally acyclic relatively to 
${\cal F}'$.

For each $y'\in Z'=
Z\times_YY'$,
we have
$a_{y'}Rf'_*{\cal F}'=
(CC{\cal F}',df')_{T^*X',y'}=0$.
Since $Y'\to Y$ is \'etale at $y$,
for each $y'\in Y'\times_Yy$,
we have
$a_yRf_*{\cal F}=
a_{y'}Rf'_*{\cal F}'$
and 
$(CC{\cal F},df)_{T^*X,X_y}=
(CC{\cal F}',df')_{T^*X',y'}$.
Thus,
by applying (\ref{eqsum}) to
$f'\colon X'\to Y'$ and ${\cal F}'$,
we obtain
$$-[Y':Y]\cdot
a_yRf_*{\cal F}=
[Y':Y]\cdot
(CC{\cal F},df)_{T^*X,y}$$
and hence (\ref{eqcf}).
\qed}

\begin{cor}[{\rm cf.\ \cite[Conjecture]{Bl}}]\label{corBl}
Let $f\colon X\to Y$
be a {\em projective} flat morphism
of smooth schemes
over a perfect field $k$.
Assume that 
$\dim X=n$,
$\dim Y=1$
and that there exists
a dense open subscheme
$V\subset Y$
such that the base change
$f_V\colon X\times_YV\to V$
is smooth.
Then, for a closed point $y\in Y$,
we have
\begin{equation}
-a_yRf_*\Lambda
=
(-1)^n{c_n}^X_{X_y}(\Omega^1_{X/Y})\cap [X].
\label{eqBl}
\end{equation}
\end{cor}

\proof{
Applying Theorem \ref{thmay}
to the constant sheaf ${\cal F}=\Lambda$
and $CC\Lambda=
(-1)^n[T^*_XX]$,
we obtain
$-a_yRf_*\Lambda
=
(-1)^n
(T^*_XX,df)_{T^*X,X_y}$.
By applying Lemma \ref{lmcE} to the right hand side
and $[f^*\Omega^1_{Y/k}
\to \Omega^1_{X/k}]$,
we obtain (\ref{eqBl}).
\qed}

\cite[Conjecture]{Bl} is stated for proper morphism.
Here we need to assume $f$ to be projective
since the index formula is known only for projective schemes.

\begin{thm}\label{thmf*}
Let $f\colon X\to Y$
be a morphism
of smooth schemes
over a perfect field $k$.
Let ${\cal F}$
be a constructible complex
of $\Lambda$-modules on $X$
and $C=SS{\cal F}$ be
the singular support.
Assume that $Y$ is {\em projective},
that $f\colon X\to Y$ is {\em quasi-projective} and
is proper on the support
of ${\cal F}$ and that we have an inequality
\begin{equation}
\dim f_\circ C
\leqq \dim Y=m.
\label{eqdim}
\end{equation}
Then, we have
\begin{equation}
CCRf_*{\cal F}=
f_!CC{\cal F}
\label{eqCCf}
\end{equation}
in $Z_m(f_\circ SS{\cal F})$.
\end{thm}

\proof{
We may assume that
$k$ is algebraically closed.
Since $X$ is quasi-projective,
there exists a locally closed
immersion $i\colon X\to P$
to a projective space $P$.
By decomposing $f$
as the composition of
the immersion
$(i,f)\colon X\to P\times Y$
and the second projection
$P\times Y\to Y$,
we may assume that $f$
is projective and smooth
by Lemma \ref{lmfun}.
Set $C=f_\circ SS{\cal F}
\subset T^*Y$.
We have
$SSRf_*{\cal F}
\subset
f_\circ SS{\cal F}=C$.
By the assumption, we have
$\dim C\leqq m$.
By the index formula
\cite[Theorem 7.13]{CC} and
Theorem \ref{thmay},
the equality (\ref{eqCCf}) is proved for
$Y$ of dimension $\leqq 1$.
We show the general case
by reducing to the case $\dim Y=1$.

We take a closed immersion
of $Y$ to a projective space
$i\colon Y\to {\mathbf P}$.
We use the notations 
${\mathbf P}
\overset p\gets Q
\overset{p^\vee}\to 
{\mathbf P}^\vee$ in (\ref{eqpv})
and let $p_X\colon
X\times_{\mathbf P}Q\to X$
be the projection.
After replacing the immersion $i$
by the composition with a Veronese
embedding if necessary,
we may assume that
the restriction to
${\mathbf P}(i_\circ C)
\subset Q=
{\mathbf P}(T^*{\mathbf P})$
of the projection
$p^\vee \colon Q\to {\mathbf P}^\vee$
is generically radicial
by the assumption $\dim C\leqq m=\dim Y$
and by \cite[Corollary 3.21]{CC}.
Let $C^\vee=p^\vee_\circ p_Y^\circ C
\subset T^*{\mathbf P}^\vee$ and
let $D$ denote the image
$p^\vee({\mathbf P}(i_\circ C))
\subset {\mathbf P}^\vee$.
By Lemma \ref{lmLeg},
Lemma \ref{lmpi2}.3
and the Bertini theorem,
there exists a line
$L\subset {\mathbf P}^\vee$
satisfying the following conditions:
The immersion
$h\colon L\to{\mathbf P}^\vee$
is properly $C^\vee$-transversal.
The morphism
$h\colon L\to{\mathbf P}^\vee$ meets
$p_X^\circ SS{\cal F}$
properly.
The axis $A_L$ of $L$
meets $Y$ transversely and
$L$ meets $D$ transversely.

Since
the blow-up of ${\mathbf P}$ at 
$A_L$ is the projection
$Q\times_{{\mathbf P}^\vee}L
\to {\mathbf P}$ and since $A_L$ meets $Y$ transversely,
the blow-up $Y'$ of $Y$ at
$Y\cap A_L$ fits in 
the cartesian diagram
$$
\begin{CD}
Y@<<< Y\times_{\mathbf P}Q@<<< Y'\\
@VVV
\hspace{-12mm}
\square
\hspace{10mm}@VVV
\hspace{-16mm}
\square
\hspace{12mm}@VVV\\
{\mathbf P}@<<<Q@<<<
Q\times_{{\mathbf P}^\vee}L\\
@.@VVV
\hspace{-16mm}
\square
\hspace{12mm}@VVV\\
@.{\mathbf P}^\vee@<<<L
\end{CD}$$
and is smooth over $k$.
We consider the cartesian diagram
\begin{equation}
\begin{CD}
X@<{p_X}<<X\times_{\mathbf P}Q
@<{h_X}<<X'\\
@VfVV
\hspace{-15mm}
\square
\hspace{10mm}
 @V{\tilde f}VV
\hspace{-13mm}
\square
\hspace{8mm}
@VV{f'}V\\
Y@<{p_Y}<<Y\times_{\mathbf P}Q
@<{h_Y}<<Y'\\
@.@V{p_Y^\vee}VV
\hspace{-13mm}
\square
\hspace{8mm}
@VV{p_L}V\\
@.{\mathbf P}^\vee @<h<<L
\end{CD}
\label{eqL}
\end{equation}
of projective smooth schemes over $k$.

By the identification
${\mathbf P}(T^*{\mathbf P}\times_{\mathbf P}Y)
\subsetneqq
{\mathbf P}(T^*{\mathbf P})
=Q=
{\mathbf P}(T^*{\mathbf P}^\vee)$,
the operation $p^\vee_{Y\circ} p_Y^\circ
=p^\vee_\circ p^\circ i_\circ$
defines an injection
from the set of irreducible
closed conical subsets of dimension $m=\dim Y$ of
$T^*Y$ and 
to the set of irreducible
closed conical subsets of 
$T^*{\mathbf P}^\vee$.
Hence the equality
(\ref{eqCCf})
is equivalent to
\begin{equation}
p^\vee_{Y!} p_Y^!CCRf_*{\cal F}=
p^\vee_{Y!} p_Y^!f_!CC{\cal F}.
\label{eqYP}
\end{equation}
It suffices to compare the
coefficients of 
$C_a^\vee=p^\vee_{Y\circ} p_Y^\circ C_a
\subset T^*{\mathbf P}^\vee$
for each irreducible component
of $C=\bigcup_aC_a$ of
dimension $m=\dim Y$.

Since the restriction 
${\mathbf P}(i_\circ C)
\to D$ of the projection
$p^\vee\colon 
Q\to {\mathbf P}^\vee$
is radicial on a dense open subset,
irreducible components
${\mathbf P}(i_\circ C_a)$
of 
${\mathbf P}(i_\circ C)$
correspond uniquely
with irreducible components $D_a$ 
of $D$ defined as the images.
Since $L$ meets $D$ transversely,
each point of the intersection
$L\cap D$ is contained in
exactly one irreducible component
$D_a$ of $D$.
Hence, (\ref{eqYP}) is further equivalent to
\begin{equation}
h^!p^\vee_{Y!} p_Y^!CCRf_*{\cal F}=
h^!p^\vee_{Y!} p_Y^!f_!CC{\cal F}.
\label{eqLCCf}
\end{equation}

Let $\pi_X\colon X'\to X$
denote the composition 
$p_X\circ h_X$ of the top line
in (\ref{eqL}).
We show that
the equality (\ref{eqLCCf})
is equivalent to
\begin{equation}
CCR(p_Lf')_*\pi_X^*{\cal F}=
(p_Lf')_!CC\pi_X^*{\cal F}.
\label{eqCCf'}
\end{equation}
First, we compare the left hand sides.
By \cite[Corollary 7.12]{CC}
applied to $i_*Rf_*{\cal F}$ on
${\mathbf P}$,
the left hand side of (\ref{eqLCCf}) equals
$h^!CCRp^\vee_{Y*} p_Y^*Rf_*{\cal F}$.
Since $SSRp^\vee_{Y*} p_Y^*Rf_*{\cal F}
\subset C^\vee$
and since $h\colon L\to {\mathbf P}^\vee$
is properly $C^\vee$-transversal,
the left hand side further equals 
$CCh^*Rp^\vee_{Y*} p_Y^*Rf_*{\cal F}$
by \cite[Theorem 7.6]{CC}.
By proper base change theorem,
this is equal to the left hand side
$CCR(p_Lf')_*\pi_X^*{\cal F}$
of (\ref{eqCCf'}).

Next, we compare the right hand sides.
The right hand side of
(\ref{eqLCCf}) is equal to
$(p_L f')_!\pi_X^!CC{\cal F}$
by the projection formula
\cite[Theorem 6.2 (a)]{Ful}.
Since $h\colon L\to {\mathbf P}^\vee$
is $C^\vee$-transversal
and $C^\vee=p^\vee_{Y\circ} p_Y^\circ C
=(p^\vee_Y \tilde f)_\circ p_X^\circ 
SS{\cal F}$,
the immersion 
$h_X\colon X'\to X\times_{\mathbf P}Q$
is 
$p_X^\circ SS{\cal F}$-transversal
by Lemma \ref{lmtrbc}.2.
Further since $h\colon L\to {\mathbf P}^\vee$
meets $p_X^\circ SS{\cal F}$ properly,
the immersion 
$h_X\colon X'\to X\times_{\mathbf P}Q$
is properly
$p_X^\circ SS{\cal F}$-transversal.
Since $p_X$ is smooth, the composition
$\pi_X=p_X\circ h_X$
is properly
$SS{\cal F}$-transversal.
Thus by \cite[Theorem 7.6]{CC},
it further equals
to the right hand side
$(p_Lf')_!CC\pi_X^*{\cal F}$
of (\ref{eqCCf'}).

We show the equality (\ref{eqCCf'})
by applying Theorem \ref{thmay}
to complete the proof.
The largest
open subset where the projection
$p^\vee_Y\colon
Y\times_{\mathbf P}Q
\to {\mathbf P}^\vee$
is $p_Y^\circ C$-transversal
equals 
the largest
open subset where the immersion
$Y\times_{\mathbf P}Q
\to Y\times {\mathbf P}^\vee$
is $C\times T^*{\mathbf P}^\vee$-transversal
and hence is the complement of 
${\mathbf P}(i_\circ C)
\subset
Y\times_{\mathbf P}Q
={\mathbf P}(Y\times_{\mathbf P}T^*{\mathbf P})$.
Since $p_Y^\circ C=
p_Y^\circ f_\circ SS{\cal F}=
\tilde f_\circ p_X^\circ SS{\cal F}
=
\tilde f_\circ SSp_X^*{\cal F}$,
the composition 
$p_Y^\vee\tilde f\colon
X\times_{\mathbf P}Q
\to {\mathbf P}^\vee$
is $SSp_X^*{\cal F}$-transversal 
on the complement $ 
{\mathbf P}^\vee\sm D$
by \cite[Lemma 3.8 (1)$\Rightarrow$(2)]{CC}.
By Lemma \ref{lmtrbc},
the morphism $p_Lf'\colon X'\to L$
is $SS\pi_X^*{\cal F}$-transversal
on the dense open subset
$L\sm L\cap D$.
Hence the equality (\ref{eqCCf'})
follows from Theorem \ref{thmay}
applied to $\pi_X^*{\cal F}$
and the equality (\ref{eqCCf})
is proved.
\qed}

\medskip

To prove Theorem \ref{thmf*},
we need the assumption that $Y$
is projective,
because we assume in Theorem \ref{thmay}
that $f$ loc.\ cit.\ is projective
on the support of ${\cal F}$.
If we can weaken this assumption on $f$
in Theorem \ref{thmay},
we will be able to remove the assumption on $Y$
in Theorem \ref{thmf*}
since the characteristic cycle is characterized
by the Milnor formula (\ref{eqMil}).


In the case of characteristic $0$,
we recover the classical result
as in \cite[Proposition 9.4.2]{KSc},
under the extra assumption that 
$Y$ is projective and that $f$ is quasi-projective.
Let $X$ be a smooth scheme
equidimensional of dimension $n$ over a field $k$
and let $\omega_X\in \Omega^2(T^*X)$ denote
the canonical symplectic form 
on the cotangent bundle $T^*X$.
Let $C\subset T^*X$ be a closed conical subset.
We say that $C$ is
{\em isotropic} if the restriction of
$\omega_X$ on $C$ is $0$.
We say that $C$ is
{\em Lagrangian} if it is isotropic and
if $C$ is equidimensional of dimension $n$.

\begin{lm}\label{lmch0}
Let $k$ be a field of characteristic $0$
and let $f\colon X\to Y$ be a morphism
of smooth schemes over $k$.
Assume that $X$ (resp.\ $Y$) is 
equidimensional of dimension $n$ 
(resp.\ $m$).
Let $C\subset T^*X$ be a closed conical subset.
If $C\subset T^*X$ is isotropic,
then $f_\circ C\subset T^*Y$ is also isotropic.
\end{lm}

The author learned the following proof
from Beilinson.

\proof{
Let $T^*_\Gamma(X\times Y)
\subset T^*(X\times Y)$
be the normal bundle
of the graph $\Gamma\subset X\times Y$
 of $f\colon X\to Y$
and let $p_2\colon T^*(X\times Y)=
T^*X\times T^*Y\to T^*Y$
be the projection.
The direct image $f_\circ C\subset T^*Y$
equals the image by $p_2$ of the intersection
$C_1=T^*_\Gamma(X\times Y)
\cap (C\times T^*Y)$.

Since the normal bundle
$T^*_\Gamma(X\times Y)
\subset T^*(X\times Y)$
is isotropic
and since
$\omega_{X\times Y}$ equals
the sum $p_1^*\omega_X+p_2^*\omega_Y$
of the pull-backs by projections,
the assumption that
$C\subset T^*X$ is isotropic implies
that the restriction of
$p_2^*\omega_Y$ on $C_1$ is 0.
Since $k$ is of characteristic $0$,
for each irreducible component
$C'$ of $f_\circ C\subset T^*Y$,
there exists a closed subset
$C'_1\subset C_1$
generically \'etale over $C'$.
Hence the assertion follows.
\qed}

\begin{pr}\label{prch0}
Let $k$ be a field of characteristic $0$
and let $X$ be a smooth scheme over $k$.
Let ${\cal F}$ be a constructible
complex of $\Lambda$-modules on $X$.

{\rm 1.}
The singular support
$SS{\cal F}$ is Lagrangian.

{\rm 2.}
Let $f\colon X\to Y$ be a morphism
of smooth schemes over $k$.
Assume that $f$ is proper on
the support of ${\cal F}$.
Then, the inequality
{\rm (\ref{eqdim})} holds.
Further if $f\colon X\to Y$ is quasi-projective,
the equality {\rm (\ref{eqCCf})} holds.
\end{pr}

\proof{
1.
We may assume that $X$
is equidimensional of dimension $n$.
Since the singular support $SS{\cal F}$ 
is equidimensional of dimension $n$
\cite[Theorem 1.3 (ii)]{Be},
it suffices to show that
$SS{\cal F}$ 
is isotropic.
By devissage,
we may assume that
there exist a locally closed immersion
$i\colon V\to X$ of smooth scheme,
a locally constant sheaf 
${\cal G}$ on $V$
and 
${\cal F}=i_!{\cal G}$.

Since the resolution of singularity
is known in characteristic $0$,
the immersion $i$ is decomposed
by an open immersion
$j\colon V\to W$ and a
proper morphism $h\colon W\to X$
such that $W$ is smooth
and $V$ is the complement of
a divisor with simple normal crossings.
Thus, by the inclusion $SS{\cal F}
=SS Rh_*j_!{\cal G}
\subset h_\circ SS j_!{\cal G}$
and Lemma \ref{lmch0},
it is reduced to the case
where $i=j$ is an open immersion
of the complement of
a divisor with simple normal crossings.
Since $k$ is of characteristic $0$,
this case is proved in
\cite[Proposition 4.11]{CC}.

2.
By 1 and Lemma \ref{lmch0},
the direct image
$f_\circ SS{\cal F}$ is
isotropic. Hence
the inequality 
$\dim f_\circ SS{\cal F}\leqq
\dim Y$ (\ref{eqdim}) holds.

We show the equality
$CCRf_*{\cal F}
=f_!CC{\cal F}$ (\ref{eqCCf}).
Similarly as in the proof of Theorem \ref{thmf*},
we may assume that $Y$ is affine
and $f\colon X=P\times Y\to Y$
is the projection for a projective smooth scheme
$P$ over $k$.
By resolution of singularity,
we may assume that $Y$ is projective and smooth.
Then since
the inequality (\ref{eqdim}) holds, 
we may apply Theorem \ref{thmf*}.
\qed}

\subsection{Index formula for vanishing cycles}\label{ssiv}

We prepare some notation
to formulate an index formula
for vanishing cycle complex.
Let $f\colon X\to Y$
be a smooth morphism
of smooth schemes
over a perfect field $k$.
Assume that $X$ 
(resp.\ $Y$) is
equidimensional of
dimension $n+1$ (resp.\ $1$).
Let ${\cal F}$ be a constructible
complex of $\Lambda$-modules
on $X$.
Let $y\in Y$ be a closed point
and $i_y\colon X_y\to X$
be the closed immersion
of the fiber.
Assume that $f\colon X\to Y$
is properly $SS{\cal F}$-transversal
on the complement $X\sm X_y$
of the fiber $X_y=f^{-1}(y)$.
Then, the specialization
\begin{equation}
{\rm sp}_ySS{\cal F}
\subset
T^*X_y
\label{eqspSS}
\end{equation}
is defined as a closed
conical subset 
equidimensional of dimension $n$.
Further, the specialization
\begin{equation}
{\rm sp}_yCC{\cal F}
\in
Z_n({\rm sp}_ySS{\cal F})
\label{eqspCC}
\end{equation}
is defined as a cycle.

\begin{lm}\label{lmspPs}
Let $f\colon X\to Y$
be a smooth morphism
of smooth schemes over a field $k$
and assume $\dim Y=1$.
Let ${\cal F}$ be a constructible
complex of $\Lambda$-modules
on $X$ and assume
that $f\colon X\to Y$
is properly $SS{\cal F}$-transversal.
Let $y\in Y$ be a closed point.
Then,
we have
\begin{equation}
SSR\Psi_y{\cal F}
={\rm sp}_ySS{\cal F}.
\label{eqSSsp}
\end{equation}
Further if $k$ is perfect, we have
\begin{equation}
CCR\Psi_y{\cal F}
={\rm sp}_yCC{\cal F}.
\label{eqCCsp}
\end{equation}
\end{lm}

\proof{
Let $i_y\colon X_y\to X$
denote the closed immersion 
of the fiber.
Then, by the assumption that
$f\colon X\to Y$
is properly $SS{\cal F}$-transversal,
we have
${\rm sp}_ySS{\cal F}
=i_y^\circ SS{\cal F}$
and
${\rm sp}_yCC{\cal F}
=i_y^! CC{\cal F}$.
Recall that the definitions
of ${\rm sp}_y$ and $i_y^!$
both involve the minus sign.

Since $f\colon X\to Y$ is 
locally acyclic relatively to
${\cal F}$ by Lemma \ref{lmtr}.2,
the canonical morphism
$i_y^*{\cal F}
\to R\Psi_y{\cal F}$
is an isomorphism.
Hence the equalities
(\ref{eqSSsp}) and (\ref{eqCCsp}) 
follow from Lemma \ref{lmh}
and \cite[Theorem 7.6]{CC}
respectively.
\qed}
\medskip

The following example
shows that the inclusion
$SSR\Psi{\cal F}
\subset {\rm sp}_ySS{\cal F}$
does not hold in general.

\begin{ex}\label{ex2}
{\rm
Let $k$ be a field of characteristic $p>2$.
Let $X={\mathbf A}^1\times {\mathbf P}^1$
and $j\colon U={\mathbf A}^1\times {\mathbf A}^1
={\rm Spec}\ k[x,y]\to X$
be the open immersion.
Let ${\cal G}$ be the locally constant sheaf
of $\Lambda$-modules of
rank $1$ on $U$ 
defined by the Artin-Schreier covering
$t^p-t=x^py^2$ and by a non-trivial character
${\mathbf F}_p\to \Lambda^\times$.
Then, the nearby cycles complex
$R\Psi_\infty{\cal F}$
is acyclic except at
the closed point $(0,\infty)$
or at degree $1$
and 
$\dim R^1\Psi{\cal F}_{(0,\infty)}
=1.$
Hence, the singular support
$SSR\Psi_\infty{\cal F}$
equals the fiber
$T^*_{(0,\infty)}X_\infty$
at the closed point
and is not a subset of
the zero-section
${\rm sp}_\infty SS{\cal F}=
T^*_{X_\infty}X_\infty$.
}
\end{ex}

Let $Z\subset X_y$
be a closed subset.
Assume that $f\colon X\to Y$
is properly $SS{\cal F}$-transversal
on the complement of $Z$.
Then,
on the complement
$X_y\sm Z$,
we have
${\rm sp}_yCC{\cal F}
=
i_y^!CC{\cal F}
=
CCi_y^*{\cal F}$
by Lemma \ref{lmsp}
and the compatibility with
the pull-back \cite[Theorem 7.6]{CC}.
Thus, the difference
\begin{equation}
\delta_y
CC{\cal F}
=
{\rm sp}_yCC{\cal F}
-
CCi_y^*{\cal F}
\label{eqdel}
\end{equation}
is defined as a cycle in
$Z_n\bigl(Z\times_X(
{\rm sp}_ySS{\cal F}
\cup
SSi_y^*{\cal F})\bigr)$
supported on $Z$.
If $Z$ is proper over $Y$,
the intersection number 
$(\delta_y SS{\cal F},T^*_{X_y}X_y)_
{T^*X_y}$ is
defined.

\begin{pr}\label{prvan}
Let $f\colon X\to Y$
be a smooth morphism of smooth schemes
over a perfect field $k$.
Assume that $X$ 
(resp.\ $Y$) is
equidimensional of
dimension $n+1$ (resp.\ $1$).
Let ${\cal F}$ be a constructible
complex of $\Lambda$-modules
on $X$.
Let $y\in Y$ be a closed point
and let $Z\subset X_y$
be a closed subset.
Assume that $f\colon X\to Y$
is properly $SS{\cal F}$-transversal
on the complement of $Z$ 
and that either of the following 
conditions {\rm (1)} 
and {\rm (2)} is satisfied:

{\rm (1)} 
$f\colon X\to Y$ is projective.

{\rm (2)}
$\dim Z=0$.

\noindent
Then,
for the vanishing cycles complex
$R\Phi_y{\cal F}$, we have
\begin{equation}
\chi(Z_{\bar k},
R\Phi_y{\cal F})
=
(\delta_y CC{\cal F},T^*_{X_y}X_y)_
{T^*X_y}.
\label{eqvan}
\end{equation}
\end{pr}

\proof{
We may assume that
$k$ is algebraically closed.

We show the case (1).
Let $v\in Y$ be a closed point
different from $y$
and let $i_v\colon X_v\to X$
be the closed immersion.
Then, since the projective
morphism $f\colon X\to Y$
is locally acyclic relative to ${\cal F}$
outside $Z$ by Lemma \ref{lmtr}.2,
the left hand side of
(\ref{eqvan})
equals 
\begin{equation}
\chi(Z,
R\Phi_y{\cal F})
=
\chi(X_y,
R\Psi_y{\cal F})
-
\chi(X_y,
i_y^*{\cal F})
=
\chi(X_v,
i_v^*{\cal F})
-
\chi(X_y,
i_y^*{\cal F})
\label{eqchv}
\end{equation}
The right hand side of (\ref{eqvan})
\begin{equation*}
(\delta_y CC{\cal F},T^*_{X_y}X_y)_
{T^*X_y}
=
({\rm sp}_y CC{\cal F},T^*_{X_y}X_y)_
{T^*X_y}
-
(CCi_y^*{\cal F},
T^*_{X_y}X_y)_
{T^*X_y}
\end{equation*}
equals
\begin{equation}
(i_v^! CC{\cal F},T^*_{X_v}X_v)_
{T^*X_v}
-
(CCi_y^*{\cal F},
T^*_{X_y}X_y)_
{T^*X_y}
\label{eqdelv}
\end{equation}
by (\ref{eqspiv}).
Since $i_v\colon X_v\to X$
is properly $SS{\cal F}$-transversal
by Lemma \ref{lmtrbc},
the right hand side of (\ref{eqchv})
equals (\ref{eqdelv})
by the compatibility with
the pull-back \cite[Theorem 7.6]{CC}
and the index formula
\cite[Theorem 7.13]{CC}.
Thus the equality
(\ref{eqvan}) is proved.

We show the case (2).
Since the formation of nearby cycles
complex commutes with base change
by \cite[Proposition 3.7]{TF},
we may assume that 
the action of the inertia group $I_y$
on $R\Psi_y{\cal F}$
is trivial.
Since the vanishing cycles functor
is $t$-exact by
\cite[Corollaire 4.6]{au},
we may assume that ${\cal F}$
is a simple perverse sheaf.

First, we consider the case
${\cal F}$ is supported 
on the closed fiber $X_y$.
By the assumption that
$f\colon X\to Y$
is properly $SS{\cal F}$-transversal
on the complement of $Z$,
the morphism
$f\colon X\to Y$
is locally acyclic relatively to
${\cal F}$ on the complement of $Z$.
Thus ${\cal F}$ is supported on $Z$
and the assertion follows in this case.

We may assume that
the restriction ${\cal F}|_{X_\eta}$
on the generic fiber is non-trivial.
Then, by Proposition \ref{prS}.2,
the morphism
$f\colon X\to Y$
is locally acyclic relatively to
${\cal F}$.
Hence by Lemma \ref{lmtrZ}.2,
the morphism $f\colon X\to Y$
is properly $SS{\cal F}$-transversal
and the assertion follows
from Lemma \ref{lmspPs}.
\qed}

\medskip
In the case (2) $\dim Z=0$,
Proposition \ref{prvan}
means
$CCR\Phi_y{\cal F}=
\delta_yCC{\cal F}$.
However,
Examples \ref{ex} and \ref{ex2}
show that one cannot expect
to have
$CCR\Psi_y{\cal F}=
{\rm sp}_yCC{\cal F}$ or
equivalently
$CCR\Phi_y{\cal F}=
\delta_yCC{\cal F}$
in general.

\end{document}